\documentclass{article}
\usepackage{amssymb}
\usepackage{amsmath}
\usepackage{amsthm}
\usepackage{verbatim}
\usepackage{color}
\usepackage{url}
\usepackage{graphicx}
\usepackage[square]{natbib}

\def\un{{\mathrm{1~\hspace{-1.4ex}l}}}

\newcommand{\BigO}{\mathcal{O}}

\newcommand{\jvRe}{\operatorname{Re}}
\newcommand{\jvIm}{\operatorname{Im}}
\newcommand{\jvexp}{\operatorname{exp}}
\newcommand{\frakp}{\mathfrak{p}}
\newcommand{\frakq}{\mathfrak{q}}
\newcommand{\frakU}{\mathfrak{U}}
\newcommand{\fraka}{\mathfrak{a}}

\newcommand{\opnm}{\operatorname}
\newcommand{\eps}{\varepsilon}
\newcommand{\kap}{\varkappa}
\newcommand{\Bff}{\mathbf}

\newcommand{\Op}[2]{\operatorname{Op}^w_{#1, #2}}

\newtheorem{theorem}{Theorem}[section]
\newtheorem{proposition}[theorem]{Proposition}

\newtheorem{lemma}[theorem]{Lemma}

\newtheorem{remark}[theorem]{Remark}

\numberwithin{equation}{section}

\author{Joe Viola}
\title{Non-Elliptic Quadratic Forms and Semiclassical Estimates for Non-Selfadjoint Operators}

\begin{document}

\maketitle

\begin{abstract}
We consider a class of pseudodifferential operators with a doubly characteristic point, where the quadratic part of the symbol fails to be elliptic but obeys an averaging assumption.  Under suitable additional assumptions, semiclassical resolvent estimates are established, where the modulus of the spectral parameter is allowed to grow slightly more rapidly than the semiclassical parameter.
\end{abstract}

\maketitle

\section{Introduction and Statement of Results}

\subsection{Quadratic forms and singular spaces}
Recently, there has been a renewed interest in the analysis of spectra and resolvents of non-selfadjoint operators with double characteristics. The study of
pseudodifferential operators with double characteristics has a long and distinguished tradition in the analysis of partial differential
operators \citep{Ho75}, \citep{HoALPDO3}, \citep{Sj74}. Recently, the point of view of semiclassical analysis, with important motivations
coming from the study of pseudospectra for nonselfadjoint operators, has produced a considerable body of work.

The simplest examples of pseudodifferential operators with double characteristics are quadratic differential operators
$$
	Q = Q(x, D_x) = \sum_{|\alpha| + |\beta| = 2} q_{\alpha\beta}x^\alpha D_x^\beta
$$
for $q_{\alpha\beta}\in\Bbb{C},$ $D_{x_j} = \frac{1}{i} \frac{\partial}{\partial x_j},$ and $\alpha,\beta$ multiindices. The spectrum of these operators in the 
elliptic case has been understood for some time \citep{Sj74}, but recent work \citep{Davies04}, \citep{Boulton02} showed that the operator norm of 
the resolvent $(Q - z)^{-1}$ for $z \in \Bbb{C}$ may exhibit rapid growth even far from the spectrum, when $z$ is taken along rays inside the range of the symbol,
$\{Q(x,\xi)\::\: (x,\xi) \in \Bbb{R}^{2d}\}.$ This is in sharp contrast to the case of any selfadjoint operator $A$.
This rapid resolvent growth was shown to be characteristic of many non-selfadjoint pseudodifferential operators \citep{DeSjZw04} and studied 
in the case of semiclassical non-selfadjoint elliptic quadratic operators \citep{PSDuke} to demonstrate that the resolvent of
$Q(x,hD_x)$ grows exponentially quickly, as the semiclassical parameter $h\rightarrow 0$, when the spectral parameter $z$ lies inside the range of 
the symbol. Rescaling shows that growth in $|z|$ along rays inside the range of $Q(x,\xi)$ when $h = 1$ is fixed is equivalent to growth in $h^{-1}$ for $z$ fixed inside the range of
$Q(x,\xi)$.

The region in $\Bbb{C}$ where the resolvent of an operator grows large is called the pseudospectrum, and the breadth of the
pseudospectrum corresponds to instability of the spectrum under small perturbations.  A natural question is to what extent
lower-order terms in the symbol of an operator with double characteristics may perturb the resolvent growth and spectrum governed
by the quadratic part.  The study of the pseudospectrum of a variety of operators has received much recent interest in a diverse
array of problems, and an overview may be found in \citep{TrefEmb05}.

Of particular relevance here are the investigations \citep{H-PS0}, \citep{H-PS} dealing with non-selfadjoint operators with double
characteristics which have non-elliptic quadra\-tic approximations at the double characteristics. Those works were directly inspired by the
fundamental paper \citep{HeSjSt05}, devoted principally to operators of Kra\-mers-Fokker-Planck type. Perhaps the simplest example of such an
operator is obtained when considering the Weyl quantization of the quadratic form
\begin{equation}\label{eKFP}
q(x,y,\xi,\eta) = \frac{1}{2}\left(y^2 + \eta^2\right) + i \left(y\xi - ax\eta\right), \quad a \in \Bbb{R}\backslash \{0\},
\end{equation}
and $(x,y,\xi,\eta)\in \Bbb{R}^4$. Notice that the quadratic form $q$ is not elliptic on $\Bbb{R}^4,$ since $q(x,0,\xi,0) = 0.$
Nevertheless, it is known that the corresponding operator $q^w$, given as the Weyl quantization of $q$, has discrete spectrum and the 
associated semigroup is well-behaved \citep{HeSjSt05}.

In the work \citep{H-PS0}, it was realized that the particular example given by quantizing (\ref{eKFP}) could be understood by means of the so-called singular space
$S$, contained in the phase space $T^*\Bbb{R}^d \sim \Bbb{R}^d_x \times \Bbb{R}^d_\xi$, and intrinsically associated to a general quadratic form $q$ on the phase space
for which $\jvRe q\geq 0$. Since the singular space $S$ will play a crucial role in the present work, we shall now pause to recall its definition,
following \citep{H-PS0}.

Recall the standard symplectic inner product on $\Bbb{R}^{2d}$ (or $\Bbb{C}^{2d}$),
\begin{equation}
\label{eqnsigmadef}
\sigma((x,\xi),(y,\eta)) = \langle \xi, y\rangle - \langle x, \eta\rangle =
\langle (x,\xi), J(y, \eta)\rangle,
\end{equation} with
\begin{equation}\label{eJDef}
J = \left(\begin{array}{cc} 0 & -I \\ I & 0
\end{array}\right).
\end{equation}
Here and throughout the paper, inner products on $\Bbb{R}^{2d}$ or $\Bbb{C}^{2d}$ will be symmetric instead of Hermitian,
meaning $\langle x, y\rangle = \sum_{j=1}^d x_j y_j$ without  taking complex conjugates.  For brevity of notation, we will
frequently use capital letters for elements of $\Bbb{R}^{2d}\sim T^*\Bbb{R}^d,$ as in \[X = (x,\xi) \in
\Bbb{R}^{2d}.\] The Hamilton map $F$ of a complex-valued quadratic form $q(X)$ is the unique complex linear $F$ for which
\[
	\sigma(X, FY) = -\sigma(FX, Y), \quad \forall X,Y \in \Bbb{R}^{2d}
\]
and for which 
\[
	q(X) = \sigma(X, FX), \quad \forall X \in \Bbb{R}^{2d}.
\]
Writing 
\begin{equation}\label{eqninnerprodpolar}
	q(X) = \langle X, AX\rangle
\end{equation}
where $A$ is symmetric gives, by (\ref{eqnsigmadef}), \[JF = A.\]  We also note here that, recalling the Hamilton vector field of a $C^1$ function $f$,
\begin{equation}\label{eDefH}
H_f = \sum_{j=1}^d \frac{\partial f}{\partial\xi_j}\frac{\partial}{\partial x_j} - \frac{\partial f}{\partial x_j}\frac{\partial}{\partial \xi_j},
\end{equation}
it is easy to check that
\begin{equation}\label{eHq2F}
H_q = 2 F
\end{equation}
for $q$ a quadratic form. The formula continues to hold if $q$ and $F$ are replaced by their respective real or imaginary parts.

The definition of the singular space $S$ of the quadratic part $q(X)$, given in \citep{H-PS0}, is
\[
S := \left(\bigcap_{k = 0}^\infty \operatorname{ker}\left[\jvRe F\circ(\jvIm F)^k\right]\right)\cap \Bbb{R}^{2d}.
\]
By the Cayley-Hamilton theorem, it is sufficient to take $k = 0, 1, \dots, 2d-1$ in the definition of $S$.  
Furthermore, following arguments in \citep{H-PS0}, we see that the singular space $S$ can be characterized as follows,
\begin{equation}
\label{eqS}
S = \left\{ X\in \Bbb{R}^{2d}; H_{\jvIm q}^k \jvRe q(X)=0,\quad k=0,1,2,\ldots\right\}.
\end{equation}
It follows that the singular space $S$ plays a natural role when
investigating whether, when $q$ is a quadratic form with positive semidefinite real part, $\jvRe q$ becomes positive definite when averaged along the Hamilton flow of the imaginary part of $q$.

In the present paper, continuing the analysis of \citep{H-PS},
we shall be concerned with semiclassical resolvent estimates for operators with double characteristics, whose quadratic approximations at doubly
characteristic points satisfy certain averaging-type conditions, expressed by means of the singular space $S$.

\subsection{Statement of main result}
Let $p\in S(1)$, where the symbol class $S(1)$ is defined as follows,
$$
S(1):= \{a\in C^{\infty}(\Bbb{R}^{2d},\Bbb{C}) \::\: |\partial^\alpha a(x,\xi)| \leq \BigO_\alpha(1)\}.
$$
We shall assume that
\begin{equation}
\label{eqnCharSetIsZero}
\jvRe p(x,\xi) \geq 0,~~~\jvRe p(x,\xi) = 0 \Leftrightarrow (x,\xi) =
(0, 0).
\end{equation}
At the point $(0,0)$, we also assume that $\jvIm p$ vanishes to second order,
\begin{equation}\label{eqnVanishesSecondOrder}
\jvIm p(0,0) = \jvIm p'(0,0) = 0,
\end{equation}
and as a consequence we will refer to $(0,0)$ as the (unique) doubly characteristic point of $p.$
We furthermore assume ellipticity at infinity of $\jvRe p$ in the class $S(1)$, meaning that
\begin{equation}\label{eqnEllipticityAtInfinity}
\liminf_{|(x,\xi)|\rightarrow\infty}\jvRe p(x,\xi) > 0.
\end{equation}

Let us consider the Taylor expansion of $p$ at the origin $(x,\xi)=(0,0)$,
\[p(x,\xi) = q(x,\xi) + \BigO(|(x,\xi)|^3).\]
Here $q(x,\xi)$ is a quadratic form such that $\jvRe q(x,\xi)\geq 0$. In this paper, we will work under the assumption that the singular space associated with $q$ is trivial:
\begin{equation}
\label{eqnSisZero}
S = \{0\}.
\end{equation}
In this case, it was shown in \citep{H-PS} that the spectrum of the
semiclassical Weyl quantization of $q$ is discrete and is well understood as a lattice.  Writing $\lambda_j$ for the eigenvalues of $F = \frac{1}{2}H_q$ with positive imaginary part, of which there are necessarily $d$ when repeating for algebraic multiplicity, we have
\begin{equation}\label{eqnSpectrumLattice}
\operatorname{Spec}q^w(x,hD_x) = \left\{\frac{h}{i}\sum_{j=1}^d \left(1 + 2k_j\right)\lambda_j\::\:
k_j \in \Bbb{N}\cup\{0\}\right\}.\end{equation}
This description is precisely the same as in the globally elliptic case -- see Theorem 3.5 in \citep{Sj74}.

Associated to the symbol $p$ is the operator $p^w(x,hD_x)$, obtained as the semiclassical Weyl quantization of $p$,
\[p^w(x,hD_x)u(x) = (2\pi
h)^{-d}\iint_{\Bbb{R}^{2d}}e^{\frac{i}{h}(x-y)\cdot\xi}p(\frac{x+y}{2}, \xi)u(y)\,dy\,d\xi.\]
Here $0<h \leq 1$ is the semiclassical parameter.

The purpose of this work is to establish the following result.

\begin{theorem}\label{thmMainTheorem}
Let $p \in S(1)$ be a symbol with doubly characteristic point at $(0,0) \in \Bbb{R}^{2d}$ and elliptic real part elsewhere, as in \textrm{
(\ref{eqnCharSetIsZero})}, \textrm{(\ref{eqnVanishesSecondOrder})}, and \textrm{(\ref{eqnEllipticityAtInfinity})}. Furthermore, assume the singular
space of $q$, the quadratic part of $p$ at $(0,0)$, is trivial, as in \textrm{(\ref{eqnSisZero})}.  Let
\[
F(h) := \frac{1}{C_0}\left(\log\log\frac{1}{h}\right)^{1/d}.
\]
We will assume
that the spectral parameter $z \in \Bbb{C}$ obeys 
\[|z| \leq hF(h)\] 
and 
\[\operatorname{dist}(z, \operatorname{Spec}(q^w(x,hD_x))) \geq h e^{-F(h)/C_1}.\]
Then, for any $\rho > 0,$ there exist $h_0 > 0$ sufficiently small and $C_0,C_1 > 0$ sufficiently large where, for $z$ as above, the resolvent
\[
	(p^w(x,hD_x)-z)^{-1}:L^2(\Bbb{R}^d)\rightarrow L^2(\Bbb{R}^d)
\]
exists and satisfies, 
\[
||(p^w(x,hD_x)-z)^{-1}||_{L^2\rightarrow L^2} \leq \BigO(h^{-1-\rho})
\]
for all $0 < h \leq h_0.$
\end{theorem}

\begin{figure}\label{f1}
\caption{Example of valid $z$ for theorem \ref{thmMainTheorem}}
\begin{center}
\includegraphics[scale=0.35]{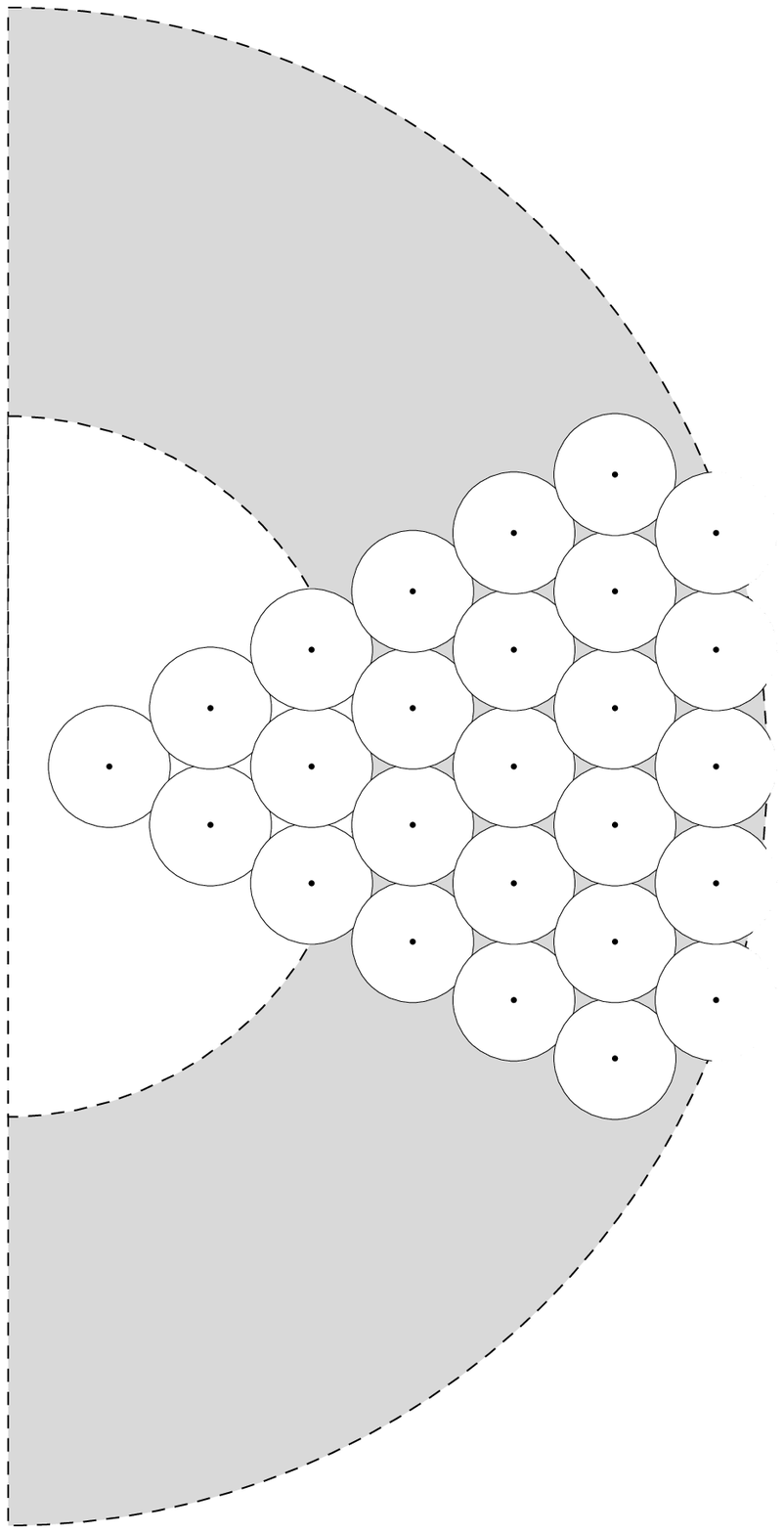}
\includegraphics[scale=0.35]{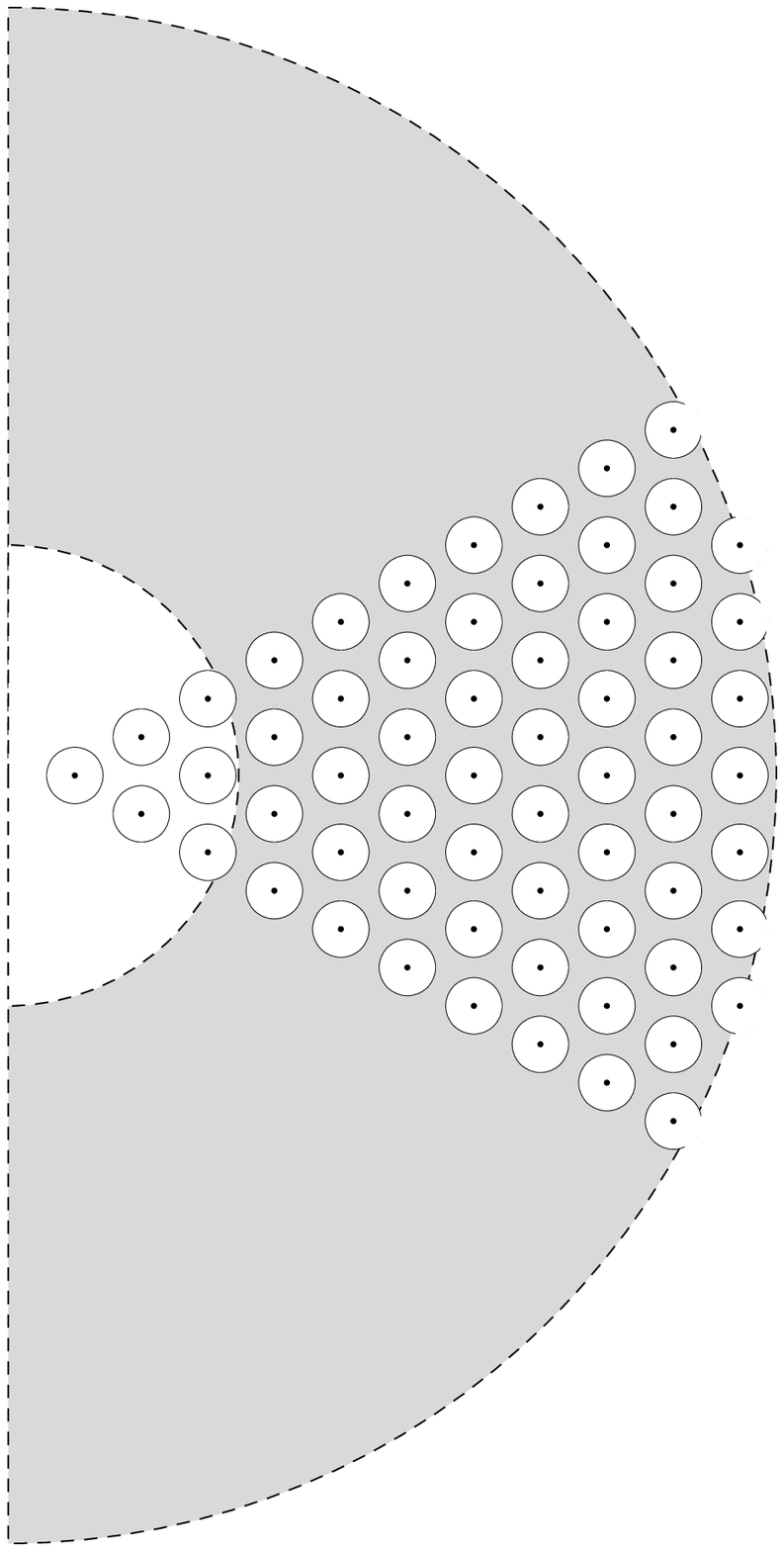}
\includegraphics[scale=0.35]{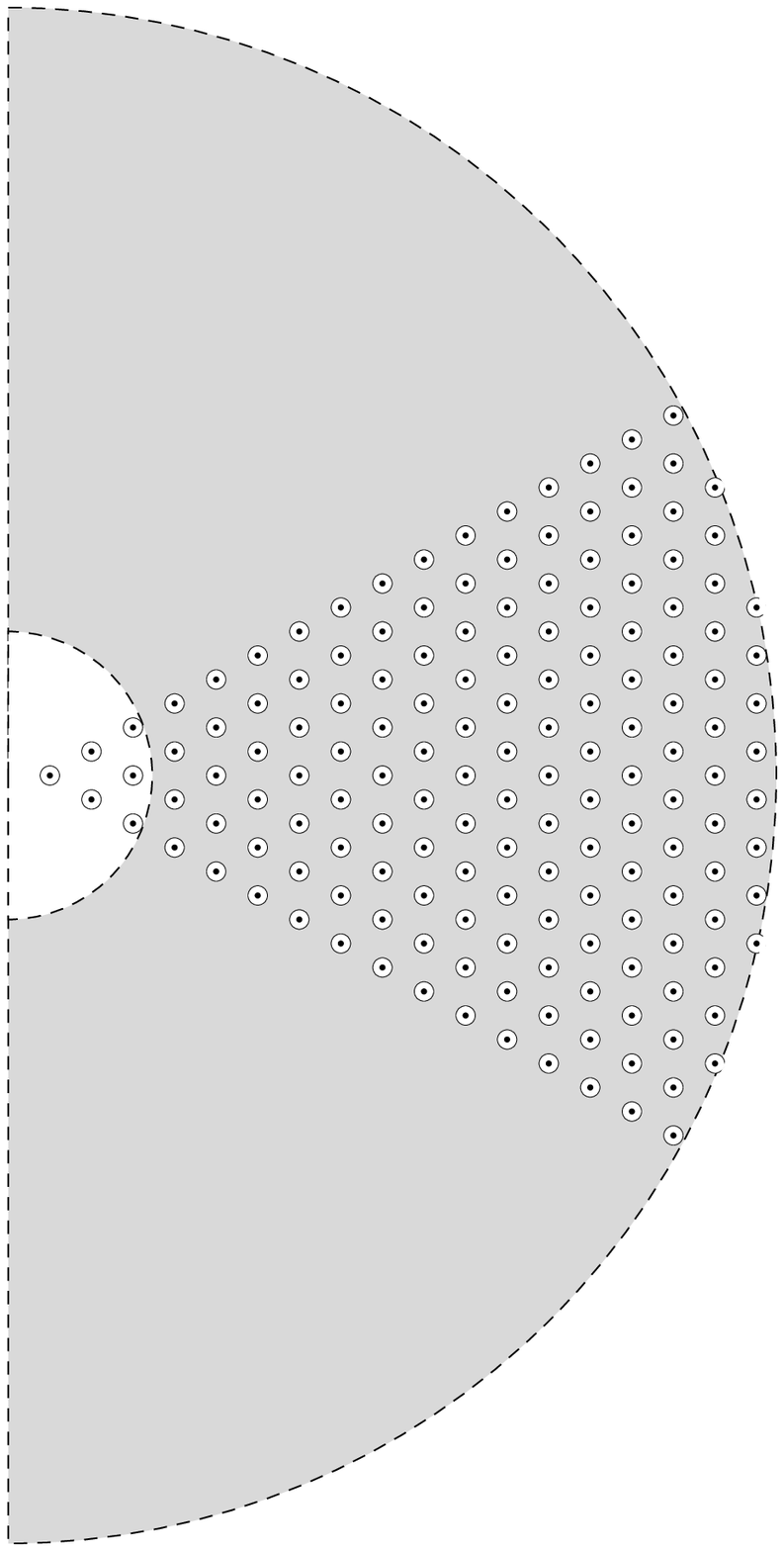}
\end{center}
\end{figure}

\begin{remark} The main novelty of theorem \ref{thmMainTheorem} is that here a polynomial resolvent bound is shown to hold when the spectral parameter $z$ may become 
$\gg h$ in modulus, although, unfortunately, it should be bounded by $hF(h)$. In the region $z=\mathcal{O}(h)$, the situation is much more pleasant, as 
was established recently in \citep{H-PS}. Specifically, assuming that the quadratic form $q$ is elliptic along the singular space $S$, it was shown 
in \citep{H-PS} that one has a global semiclassical resolvent estimate,
\begin{equation}
\label{res}
\left(p^w(x,hD_x)-h z\right)^{-1}  = \mathcal{O}\left(h^{-1}\right): L^2 \rightarrow L^2,
\end{equation}
provided that the spectral parameter $z\in \Bbb{C}$ varies in a bounded region while avoiding the spectrum of $q^w(x,D_x)$. The case when $q$ is 
globally elliptic is classical and goes back to \cite{Sj74}. The paper \citep{H-PS} extended (\ref{res}) beyond the elliptic 
case while still assuming that the spectral parameter is confined to a region of size $h$ around the origin. The restriction to such 
small $h$-dependent sets is dictated here by the scaling properties of the quadratic part $q$. Moreover, going beyond such a region, one enters 
deeper into the pseudospectrum, where spectral instability takes over and some wild resolvent growth is known, even in the purely quadratic 
case \citep{PSDuke}. It is therefore a natural and challenging problem to show that the resolvent may be polynomially controlled in regions of 
size asymptotically larger than $h,$ which is precisely the subject of this work.

In figure \ref{f1}, we have an illustration of a typical region in $\Bbb{C}$ to which the theorem applies, for decreasing values of $h$.  The grey wedge surrounded by the dashed line is $\{3h \leq |z| \leq hF(h), \jvRe z \geq 0\},$ corresponding to an operator such as Kramers-Fokker-Planck whose symbol has a range of the entire right half-plane.  The points in a lattice in the region $\{|\arg z| \leq \pi/6\}$ represent the spectrum of $q^w(x,hD_x)$ according to (\ref{eqnSpectrumLattice}), here assuming that $\opnm{Spec}F = \{\pm e^{i\pi/3}, \pm e^{2i\pi/3}\}$.  The circles surrounding those points are the forbidden region
\[
\operatorname{dist}(z, \operatorname{Spec}q_1^w(x,hD_x)) \leq he^{-F(h)/C_1},
\]
with $C_1 = 10$.  The inner white region is of order $h,$ in fact, of radius $3h,$ which we recall has already been addressed by \citep{H-PS}.  The grey region extends to $hF(h)$ for $F(h) = 6.5, 10, 16.$  Note that values of $F(h)$ correspond to extremely small $h$, of size $h = \jvexp(\jvexp(-CF(h)^d)).$   Furthermore, note that the excepted discs surrounding the spectrum occupy a vanishing fraction of the allowed grey region, which is stated precisely and proven in the body of the paper.
\end{remark}

\begin{remark} In \citep{Vi09}, the author proves a similar result in the case where the quadratic part $q$ of the symbol $p$ is globally elliptic. The corresponding main result is stronger in that it applies to $z$ in a region with a larger expansion factor, $|z| \leq hf(h)$ with
\[
	f(h) = \frac{1}{C}\left(\frac{\log \frac{1}{h}}{\log\log\frac{1}{h}}\right)^{1/d}.
\]
(We remark that the denominator could be removed using the same improvement presented in section \ref{sImprovedQuad}.)
The essential difference is that, in the elliptic case, only one rescaling is used, putting the slightly larger spectral parameter
$|z| \sim hf(h)$ into a fixed bounded region and inducing the rescaled semiclassical parameter $\tilde{h} = 1/f(h).$  In the current context,
with subelliptic quadratic part, the FBI transform methods following \citep{HeSjSt05} make the symbol elliptic only in a region of size approximately
$(h \log \frac{1}{h})^{1/2},$ at which point an additional cutoff and rescaling must be applied.  In the end, the error due to the shift of contour
on the FBI transform side, which makes the quadratic part of $p$ elliptic, dominates the error from estimating the resolvent of the new, elliptic symbol.

The methods from \citep{Vi09} are critical in providing upper bounds for the resolvent in section \ref{sLocalQuad},
once the symbol is made elliptic using the weight function in section \ref{SecWeightFunc}.
\end{remark}

\begin{remark} It is possible to extend the statement of theorem \ref{thmMainTheorem} to suitable symbols $p\in S(m)$, where $m\geq 1$ is an order function,
provided that the assumptions (\ref{eqnCharSetIsZero}) and (\ref{eqnVanishesSecondOrder}) hold, and that the ellipticity hypothesis
(\ref{eqnEllipticityAtInfinity}) is modified
accordingly. Also, the point $(0,0)$ in (\ref{eqnCharSetIsZero}) can be replaced by an arbitrary finite set $\subset \Bbb{R}^{2d}$, provided that the
singular space for the quadratic approximation at each doubly characteristic point is trivial. These remarks, made in analogy with the end of section 1 in \citep{H-PS}, can
be inferred following \citep{H-PS} as well as the methods of the present paper.
\end{remark}

\begin{remark} It seems quite likely that a major part of the following analysis goes through if we make the weaken the assumption (\ref{eqnSisZero}) to require only that the quadratic form $q$ is elliptic along $S$ in the sense that $(\jvRe q)^{-1}(\{0\}) \cap S = \{0\}$. We intend to return to this observation later,
and hope to treat this more general situation in a future paper.
\end{remark}

The plan of the paper is as follows. In section \ref{SecWeightFunc}, we establish the weight function $G_\eps$ and its quadratic part
$G_q,$ and the associated IR-submanifold of $\Bbb{C}^{2d}$ along which the symbol $p$ is elliptic. In section \ref{sFBI}, we recall the
well-known tools associated with the FBI transform upon which the proof will rely, including change of contour with small error and a localized
quantization-multiplication formula. In section \ref{sLocalQuad}, we demonstrate local resolvent estimates for the now-elliptic quadratic part
of $p,$ in a weighted space, and in section \ref{sLocalFull} we prove a result which extends these results from the quadratic part to the full symbol.
Since we are aiming at getting resolvent bounds that are polynomial, we shall take $\eps$ only logarithmically larger than $h$, 
so that when comparing the weighted and the unweighted norms, only polynomial-in-$h$ losses are obtained.
Section \ref{sExterior} provides the corresponding estimates in the exterior region, and they are then glued together to prove the main theorem 
in section \ref{sProof}.

\subsection{Examples}

For a first example, consider the operator $q^w(x, y, hD_x, hD_y)$ of Kramers-Fokker-Plank type with symbol given by (\ref{eKFP}), studied previously in \citep{HeSjSt05}.   We can easily compute that, for this $q,$
$$
F = \left(\begin{array}{cccc} 0 & -ia/2 & 0 & 0
\\ i/2 & 0 & 0 & 1/2
\\ 0 & 0 & 0 & ia/2
\\ 0 & -1/2 & -i/2 & 0 \end{array}\right),
$$
from which $\operatorname{ker}(\jvRe F) = \{(x,0,\xi,0)\}$ and $\operatorname{ker}(\jvRe F \circ \jvIm F) = \{(0, y, 0, \eta)\}.$  The intersection which defines $S$ is therefore $\{0\}.$

For a second example, we shall follow \citep{H-PS} and construct a symbol in $S(1)$ when $d = 2.$ Let $V$ and $W$ be two $C^{\infty}_b(\Bbb{R}^2,\Bbb{R})$ functions.  We make the assumption that the non-negative function $V \geq 0$ is elliptic at
infinity in the sense of (\ref{eqnEllipticityAtInfinity})
and vanishes only when $x=0$. Focusing on the expansion of $V$ and $W$ near zero, we furthermore assume that
$$V(x)=x_1^2+\BigO(x^3)$$
and
$$W(x)=\alpha x_1^2 + 2\beta x_1 x_2 +\gamma x_2^2+\BigO(x^3),$$
when $x\rightarrow 0$, for some constants $\alpha$, $\beta$, $\gamma\in \Bbb{R}$, not all equal to zero.

To create a bounded symbol, let $\chi(\xi):\Bbb{R}^2 \rightarrow \Bbb{R}$ be a $C_0^\infty$ cutoff function taking values in $[0,1]$ and equal to 1 in a neighborhood of 0.  Considering the symbol
$$p(x,\xi)=\chi(\xi)\xi^2 + (1-\chi(\xi))+V(x)+iW(x),$$
we notice that
$$(\jvRe p)^{-1}(0)=\{(0,0,0,0)\},$$
and that this symbol satisfies the assumptions (\ref{eqnCharSetIsZero}), (\ref{eqnVanishesSecondOrder}), and (\ref{eqnEllipticityAtInfinity}).
The quadratic approximation of $p$ at $(0,0,0,0)$ is then given by the following quad\-ra\-tic form
\begin{equation}\label{ex1}
q(x_1,x_2,\xi_1,\xi_2)=\xi_1^2 + \xi_2^2+ x_1^2 + i(\alpha x_1^2 + 2\beta x_1 x_2 +\gamma x_2^2).
\end{equation}
Precisely when $\gamma \neq 0,$ we have $q(x_1, x_2, \xi_1, \xi_2) = 0$ only at $0 \in \Bbb{R}^4.$  In view of Lemma 3.1 in \citep{Sj74}, this suffices to show $q$ is elliptic in the additional sense that $q(\Bbb{R}^4)$ is a closed proper cone in $\Bbb{C},$ used in \citep{Sj74}, \citep{HeSjSt05}, \citep{H-PS0}, \citep{H-PS}, \citep{Vi09}.  When $\gamma = 0,$ we see that $\operatorname{ker}\jvRe F = \{(0, x_2, 0, 0)\},$ and that $\operatorname{ker}(\jvRe F \cdot \jvIm F) = \{(0, 0, \xi_1, \xi_2)\}$ if $\beta \neq 0$ and contains $(0, x_2, 0, 0)$ when $\beta = 0.$  Since $\jvRe F (\jvIm F)^2 = 0,$ we conclude that $S = \{0\}$ only when either $\gamma = 0$ or when $\gamma = 0, \beta \neq 0.$ In these cases, theorem \ref{thmMainTheorem} can therefore be applied to the operator $p^w(x,hD_x)$

In the case when $\beta=\gamma=0$, the singular space
$S = \{(0, x_2, 0, 0)\}.$  In this case for our example $p,$ the singular space $S$ is precisely equal to the kernel of the full
Hamilton map $F$ for $q,$ and so theorem \ref{thmMainTheorem} does not apply here. In fact, the spectrum of the associated operator
$$q^w(x,D_x)=D_{x_1}^2 + D_{x_2}^2+ (1+i\alpha)x_1^2,$$
is no longer discrete.

\medskip
\noindent\textbf{Acknowledgement:} The author would like to express his gratitude to Michael Hitrik for his guidance and help in preparing this work.

\section{Exponential weights on the phase space}\label{SecWeightFunc}

\subsection{The quadratic case}\label{SecWeightFuncQuadratic}

We recall the standard definition of the Hamilton vector field $H_f$ for $f \in C^1(\Bbb{R}^{2d}, \Bbb{C})$, given by (\ref{eDefH}).  In the case where $q$ is quadratic, we also recall the relation (\ref{eHq2F}) between the Hamilton vector field of $q$ and the matrix $F$, and that the relation remains true after taking real and imaginary parts.

If $f, g: \Bbb{R}^{2d}\rightarrow \Bbb{R}$ are sufficiently regular functions, we define the average of $f$ over the Hamilton flow of $g$ up to
time $T\neq 0$ via \[\langle f\rangle_{g, T}(X) := \frac{1}{T}\int_0^T f(\jvexp (tH_g)X)\,dt.\]  The Hamilton vector field $H_g$ is defined in (\ref{eDefH}).

We will apply this average now to $\langle \jvRe q\rangle_{\jvIm q, T}.$  Since the coefficients of the vector field here are linear in $(x,\xi)$, we can conveniently alternate between viewing $\jvexp tH_{\jvIm q}$ as a solution to an ODE and as a Taylor expansion for the exponential of a matrix.  The relationship between the singular space $S$ and the positivity of averages of $\jvRe q$ along $H_{\jvIm q}$ is made precise by the following lemma.

\begin{lemma}\label{lemZeroSPositiveFlow}
Let $q:\Bbb{R}^{2d}\rightarrow \Bbb{C}$ be a quadratic form with Hamilton map $F$, obeying $\jvRe q(X) \geq 0$
for all $X \in \Bbb{R}^{2d}.$  Then the following conditions are equivalent:
\renewcommand{\labelenumi}{(\roman{enumi})}
\begin{itemize}
  \item[(i)] For any $T > 0,$
  \[\langle \jvRe q\rangle_{\jvIm q, T} > 0,\] in the sense of positive definite quadratic forms.
  \item[(ii)] $S = \{0\}.$
\end{itemize}
\end{lemma}

\begin{proof}

We remark that replacing $q$ by $\bar{q}$ changes the sign of $\jvIm F$ which preserves $S,$ from which we can see that distinctions between $T > 0$ and $T < 0$ are immaterial.  Furthermore, by (\ref{eHq2F}), we may proceed replacing $H_{\jvIm q}$ with $\jvIm F$ which suffices to prove the lemma.

The fact that $\langle \jvRe q\rangle_{\jvIm q, T}$ is a quadratic form follows immediately from linearity of $H_{\jvIm q} = 2 \jvIm F,$
that exponentials of linear maps are linear, that composition of a quadratic form with a linear map is a quadratic form, and that an
integral of a quadratic form is a quadratic form.

The fact that the second condition implies the first has been established in Proposition 2.0.1 of \citep{H-PS0}, and is seen most directly using (\ref{eqS}).
Conversely, if $S \neq \{0\},$ then there exists a nonzero $X_0 \in \Bbb{R}^{2d}$ such that $(\jvIm F)^k X_0\in
\operatorname{ker}\jvRe F$ for all $k.$  But then $e^{t\jvIm F}X_0 \in \operatorname{ker}\jvRe F$ for all $t,$
from which \[\jvRe q(e^{t\jvIm F}X_0) = \sigma(e^{t\jvIm F}X_0, (\jvRe F) e^{t\jvIm F}X_0) = \sigma(e^{t\jvIm
F}X_0, 0) = 0.\]  This shows that equivalence of the conditions.
\end{proof}

\begin{remark} Since $\jvRe q$ is a non-negative quadratic form, we may write \[\jvRe q(Y) = \langle Y, AY\rangle\] for $A$ real symmetric.  Finding $\mathbf{e}_1,\dots,\mathbf{e}_{2d} \in \Bbb{R}^{2d}$ an orthonormal basis diagonalizing $A$, we have
\[
	\jvRe q(a_1\mathbf{e}_1 + \dots + a_{2d}\mathbf{e}_{2d}) = \lambda_1 a_1^2 \mathbf{e}_1 + \dots + \lambda_{2d} a_{2d}^2\mathbf{e}_{2d}
\]
for $\lambda_j \geq 0.$  Thus for $Y \in \Bbb{R}^{2d}$ we see that $\jvRe q(Y) = 0$ if and only if $AY = 0$.  But $\jvRe F = -JA$ for $J$ in (\ref{eJDef}), and since $J$ is invertible, we conclude that
\[
	\{Y \in \Bbb{R}^{2d} \::\: \jvRe q(Y) = 0\} = \{Y \in \Bbb{R}^{2d} \::\: (\jvRe F)Y = 0\}.
\]

Now assume that $S=\{0\}$ and let us take $X \in \Bbb{R}^{2d}\backslash \{0\}$. Let $k = k(X)$ be the smallest integer
for which $(\jvRe F)\circ(\jvIm F)^k X \neq 0$.  By the previous paragraph,
\[\jvRe q((\jvIm F)^{k}X) > 0.\]
We know that the function $t\mapsto \jvRe q(e^{t\jvIm F}X)$ does not vanish to an infinite order at $t=0$, and we shall compute the leading
term in the Taylor expansion at $t=0$. To that end, notice that the fact that $(\jvRe F)(\jvIm F)^j X = 0$ for $j = 0,\dots, k-1,$ bilinearity of the symplectic form $\sigma,$ and the fact that the Hamilton map $F$ of $q$ is skew-symmetric with respect to $\sigma$ together allow us to write
\begin{multline*}
\jvRe q(e^{t\jvIm F}X) = \sigma (\frac{t^{k}}{k!}(\jvIm F)^{k}X, (\jvRe
F)\frac{t^{k}}{k!}(\jvIm F)^{k}X) + \BigO(||X||^2 t^{2k+1}) \\ =
\frac{t^{2k}}{(k!)^2}\jvRe q((\jvIm F)^{k}X) +
\BigO(||X||^2 t^{2k+1}),
\end{multline*}
where the implicit constants depend only on $k, ||\jvRe F||, ||\jvIm F||,$ and an upper
bound for $|t|,$ say $|t| \leq 1.$
\end{remark}

To exploit the positivity of the time average of $\jvRe q \circ \jvexp(tH_{\jvIm F})$ in Lemma \ref{lemZeroSPositiveFlow}, we shall employ the
method of introducing an exponential weight on the phase space $\Bbb{R}^{2d}$, following \citep{HeSjSt05} and
\citep{H-PS}.  In this section we begin with a weight which is the exponential of a real quadratic function.  The advantage of the weight $G_q$ will be that, along an approximation to the complex Hamilton flow $\jvexp(i\delta H_{G_q})$,
for small $\delta,$ we gain ellipticity in the real part of $\tilde{q}(X) \approx q(\jvexp(i\delta H_{G_q})X).$

\medskip
\noindent Let $J:\Bbb{R}\rightarrow\Bbb{R}$ be the compactly supported piecewise affine function with \[J'(t) = \delta(t) - \un_{[-1,0]}(t)\]
where $\delta$ is the Dirac mass at $0$ and $\un_{[-1,0]}$ is the characteristic function of the compact interval $[-1,0]$. We define the
real-valued quadratic form \begin{equation}\label{eqnDefGq} G_q(X) = -\int_{\Bbb{R}} J(-\frac{t}{T}) \jvRe q(e^{tH_{\jvIm q}}X)\,dt,
\end{equation}
and note that
\[
H_{\jvIm q}\jvRe q(e^{tH_{\jvIm q}}X) = \frac{d}{dt} \jvRe q(e^{tH_{\jvIm q}}X).
\]
Passing the differential operator in $X$ inside the integral and integrating by parts gives
\begin{equation}\label{eqnComputeQuadraticAveraging}
H_{\jvIm q}G_q(X) = -\frac{1}{T}\int J'(-\frac{t}{T})\jvRe q(e^{tH_{\jvIm q}}X)\,dt = \langle \jvRe q\rangle_{\jvIm q,T} - \jvRe q.
\end{equation}

We are in a position to introduce an IR-deformation of the real phase space $\Bbb{R}^{2d}$, associated to the quadratic weight $G_q$. Let
\begin{equation}\label{eqnDefLambdaDelta}
\Lambda_{\delta} = \left\{X+i\delta H_{G_q}(X):\, X\in \Bbb{R}^{2d}\right\}\subset \Bbb{C}^{2d}.
\end{equation}
Here $\delta>0$ is sufficiently small but fixed.

We shall consider the restriction of the (entire holomorphic) function $q$ to $\Lambda_{\delta}$,
\begin{equation}
\label{eqhat}
\hat{q}(X) := q(X + i\delta H_{G_q} X).
\end{equation}
Taylor expanding $q(X + i\delta H_{G_q}(X))$ to first order in $\delta$, we immediately see that, modulo $\BigO(\delta^2|X|^2),$
\begin{eqnarray*}\jvRe \hat{q} (X) &=&
\jvRe q(X) + \jvRe[i\delta(\nabla q(X) \cdot H_{G_q} (X))] \\ &=& \jvRe q - \delta \jvIm H_{G_q} q(X) \\ &=& \jvRe q(X) +
\delta H_{\jvIm q}G_q(X) \\ &=& (1-\delta)\jvRe q(X) + \delta \langle \jvRe q \rangle_{\jvIm q, T}.\end{eqnarray*}
In this computation, we switched to the differential operator perspective of $H_{G_q}$ to use that
$$
\nabla q \cdot H_{G_q} = \langle (\partial_x q, \partial_\xi q), (\partial_\xi G_q, -\partial_\xi G_q)\rangle = H_{G_q} q = - H_q G_q.
$$
Since $G_q$ is real-valued, $\jvIm H_q G_q(X) = H_{\jvIm q}G_q(X),$ and, finally, we used the identity (\ref{eqnComputeQuadraticAveraging}) for $H_{\jvIm q}G_q(X).$

We conclude that, provided that $\delta>0$ is sufficiently small,
\[
\jvRe \hat{q}(X) > 0~~\forall X \in \Bbb{R}^{2d}\backslash \{0\}.
\]

\noindent\textit{Remark}. While we may easily check that $\Lambda_\delta$ is an IR-submanifold of $\Bbb{C}^{2d},$ the map
$$
K(X) = X + i\delta H_{G_q} X
$$
need not be canonical, and it is perhaps of interest to have a canonical transformation relating the two real linear symplectic vector spaces
$\Bbb{R}^{2d}$ and $\Lambda_\delta.$ Antisymmetry of $H_{G_q}$ with respect to $\sigma$ allows us to write
$$
\sigma(K(X), K(Y)) = \sigma(X,Y) - \delta^2 \sigma(H_{G_q}X, H_{G_q}Y).
$$
We consider, for the sake of argument, $G_q(x, \xi) = x^2 + \xi^2,$ and note that, since $H_{G_q} = -2J,$ in this case $\sigma(H_{G_q}X, H_{G_q}Y) = 4\sigma(X, Y),$ and therefore we see that $K(X)$ is generally not canonical.

Instead write
$$
\sigma(K(X), K(Y)) = \sigma(X, (1 + \delta^2 H_{G_q}^2)Y) =: \sigma(X, TX)
$$
where $T$ is symmetric with respect to $\sigma.$  Let
$$
S := T^{-1/2} = (1+\delta^2 H_{G_q}^2)^{-1/2}
$$
via the Taylor series for $(1+x)^{-1/2}$ near $x = 0.$  Since $H_{G_q}^2$ is symmetric with respect to $\sigma$ and commutes with $H_{G_q},$ we note that $S$ commutes with $T$ and is symmetric with respect to $\sigma,$ and by comparing uniformly convergent power series, $S^2T = 1.$  Furthermore note that $S$ is real and close to the identity for $\delta$ sufficiently small.  Letting
\begin{equation}
\label{ekapq}
\kap_q(X) = (1+i\delta H_{G_q})S X,
\end{equation}
we exploit our computation for $K$ and obtain
$$
\sigma(\kap_q(X), \kap_q(Y)) = \sigma (S X, T SY) = \sigma(X, Y).
$$
Therefore $\kap_q$ is linear and canonical, has the same range $\Lambda_\delta$ as $K$ by closeness of $S$ to the identity, and
$$
\tilde{q}(X) = q(\kap_q X)
$$
is elliptic as $\hat{q}(X)$ is.  Furthermore, the Hamilton maps of $q$ and $\tilde{q}$ are related via the similarity relation $\tilde{F} = \kap_q^{-1}F\kap_q$ and therefore their eigenvalues are identical.

\subsection{The bounded weight function in the general case}

In what follows, we shall have to work in a microlocal exponentially weighted space, associated to a suitable weight function $G_{\eps}$, $\eps>0$,
constructed in Section 2 of \citep{H-PS}. When restricting our attention to the region where $|{X}| = \BigO(\eps^{1/2})$, the function $G_{\eps}$ is very
close to the quadratic weight $G_q(X)$ defined in Section 2.1, while further away from this region, one needs essentially to flatten out the weight so that, uniformly on
$\Bbb{R}^{2d}$, one has $G_{\eps} = \BigO(\eps)$.

The precise construction of the weight $G_{\eps}$ has been given in \citep{H-PS}, and here we shall merely describe its properties in the following
proposition, established in proposition 3 of \citep{H-PS}. In the formulation of the result, we shall simplify the statement from \citep{H-PS} slightly, taking advantage
of the fact that the singular space $S=\{0\}$ in our case, as well as the fact that the doubly characteristic set is assumed to be the single point $0 \in \Bbb{R}^{2d}$.

In section \ref{ssNearWeights} we enact this microlocal weight to improve the properties of our symbol via a shift of contour on the FBI transform side.

\begin{proposition}\label{propGeps} Let $p(x,\xi)$ stand for an almost analytic
extension of the symbol $p:\Bbb{R}^{2d}\rightarrow \Bbb{C}$,
to a tubular neighborhood of $\Bbb{R}^{2d} \subseteq \Bbb{C}^{2d}$,
which satisfies $\partial^\alpha p = \BigO_\alpha(1)$ for all $\alpha$.
Assume that $p$ continues to obey the assumptions in theorem \textrm{\ref{thmMainTheorem}}. Then there exist
constants \[C > 1,~\tilde{C} > 1,~0 < \delta_0 \leq 1,~0 <
\varepsilon_0 \leq 1\] and a weight function $G_\varepsilon \in
C_0^\infty(\Bbb{R}^{2d}; \Bbb{R})$ depending on $\varepsilon \in
(0,\varepsilon_0]$ and supported in a neighborhood of $(0,0)\in \Bbb{R}^{2d}$ such
that, uniformly for $0 < \varepsilon \leq \varepsilon_0$ and $0 <
\delta \leq \delta_0,$
\begin{itemize}
\item $G_\varepsilon = \BigO(\varepsilon),$ and $\partial^2
    G_\varepsilon = \BigO(1)$ on $\Bbb{R}^{2d}$
\item $\nabla G_\varepsilon = \BigO(|X|)$ in the region $|X|
    \leq \varepsilon^{1/2}.$
\item $\nabla G_\varepsilon = \BigO(\varepsilon^{1/2})$ in the
    region where $|X| \geq \varepsilon^{1/2}.$
\item We have
\[
\jvRe \left(p(X + i\delta H_{G_\varepsilon}(X))\right) \geq
    \frac{\delta}{\tilde{C}}\operatorname{min}(|X|^2,\varepsilon)\]
    in the region $\{|X| \leq 1/C\}$.
\item We have $$
\jvRe\left(p(X + i\delta H_{G_\varepsilon}(X))\right) \geq \frac{\delta \varepsilon}{\tilde{C}}
$$ in the region
    where $|X|^2 \geq \varepsilon.$
\end{itemize}

\end{proposition}

While we refer the reader to \citep{H-PS} for the proof, we here sketch the central ideas and the definition of $G_\eps$.  In effect, one wishes to replace $q$ in (\ref{eqnDefGq}) defining $G_q$ with the full symbol $p$.  Proposition 2 of \citep{H-PS} proves the local closeness of averages using $p$ and using $q$:
\[
	\langle \jvRe p\rangle_{T, \jvIm p} = \langle \jvRe q\rangle_{T, \jvIm q} + \BigO(|X|^3).
\]
In order to attenuate the effect of $p$ outside the small neighborhood $\{|X| \leq \eps^{1/2}\}$, one replaces $\jvRe p$ with 
\[
	(\jvRe p)_\eps(X) = g(\frac{|X|^2}{\eps}) \jvRe p(X),
\]
where $g \in C^\infty(\Bbb{R}_+, [0,1])$ is a decreasing function obeying $g(t) = 1$ for all $t \in [0,1]$ and $g(t) = 1/t$ for all $t \geq 2$.  Furthermore, one chooses $T > 0$ sufficiently small but fixed independent of $\delta, \eps$.  One then defines, in analogy with (\ref{eqnDefGq}),
\[
	G_\eps(X) = -\int_{\Bbb{R}} J(-\frac{t}{T})(\jvRe p)_\eps(e^{tH_{\jvIm p}}X)\,dt.
\]
We furthermore remark that clearly $G_\eps(0) = 0$.

In the application here, we will choose $\varepsilon$ only logarithmically larger than the semiclassical parameter
$h,$ \[\varepsilon := h \cdot \frac{1}{C}\log\frac{1}{h},\] for $C>0$ to be chosen.

\section{FBI transform tools}\label{sFBI}

The FBI transform presents an isomorphism between the space of $L^2$-functions
on a $d$-dimen\-sio\-nal real space and the space of holomorphic functions on
$\Bbb{C}^d$ which also obey an integrability condition. When passing from the Weyl quantization
on the real side to the analogous quantization on the FBI transform side, one encounters symbols defined along suitable totally real linear submanifolds
of $\Bbb{C}^{2d}$, of real dimension $2d$, as well as the corresponding contour integrals. It is through shifting the contour that we improve properties of our symbols, following \citep{H-PS} and \citep{HeSjSt05}.

An introduction to the FBI transform may be found in \citep{Martinez}; see also sections 12.2-12.5 of \citep{SjLoR}.

\subsection{FBI transforms and weighted spaces $H_{\Phi}$ for quadratic
phases}\label{subsecFBIQuadraticPhi}

Recall that the FBI transform for a holomorphic quadratic form $\varphi:\Bbb{C}^{d}_x\times \Bbb{C}^d_y \rightarrow \Bbb{C}$ which obeys
$\jvIm \varphi''_{yy} > 0$ and $\operatorname{det} \varphi''_{xy} \neq 0$ is given by $$Tu(x) = C_\varphi
h^{-3d/4}\int e^{\frac{i}{h}\varphi(x,y)}u(y)\,dy,\quad C_{\varphi}>0.$$  Where emphasis on choice of $\varphi$ is desired, we will
write $T_\varphi,$ and otherwise $T$ will, by default, refer to the standard $T_{\varphi_0}$, with
\begin{equation}\label{ephi0}
\varphi_0(x,y) =\frac{i}{2}(x-y)^2.
\end{equation}

The range of the FBI transform on $L^2(\Bbb{R}^d)$ is the space of holomorphic functions on $\Bbb{C}^d$ which are square integrable with
respect to a certain weight.
For $x \in \Bbb{C}^d,$ define
\begin{equation}\label{eqnPhiQuadraticDef}
\Phi(x) = \sup_{y\in\Bbb{R}^d}-\jvIm \varphi(x,y),
\end{equation}
which is a real-valued quadratic form. Since $\jvIm \varphi''_{yy}> 0,$ it is easy to see that $\Phi(x)$ is the unique critical value of
$y \mapsto -\jvIm \varphi(x,y).$  For $\varphi_0$ above in (\ref{ephi0}), we have the standard expression,
\begin{equation}\label{ePhi0}
\Phi_0(x) = \frac{1}{2}(\jvIm x)^2.
\end{equation}
We then define
\[
H_\Phi(\Bbb{C}^d; h) := \operatorname{Hol}(\Bbb{C}^d) \cap L^2_{\Phi}(\Bbb{C}^d;h),
\]
with $\operatorname{Hol}(\Bbb{C}^d)$ denoting entire functions on $\Bbb{C}^d$ and with the weighted space given by the norm
\[
	||v||_{L^2_\Phi(\Bbb{C}^d;h)}^2 = \int_{\Bbb{C}^d} |v(x)|^2 e^{-\frac{2}{h}\Phi(x)}\,dL(x).
\]
Here and afterwards, $dL(x)$ denotes Lebesgue measure, $d(\jvRe x)\,d(\jvIm x).$ For brevity we often omit $(\Bbb{C}^d; h)$ and write $H_\Phi$ instead.  When the semiclassical parameter needs to be emphasized (when rescaling, for example) we write $H_{\Phi, h}.$

The Weyl quantization on the FBI-Bargmann side for a quadratic weight $\varphi$ can be performed through a contour integral,
$$
\Op{\Phi}{h}(\frakp)(u)(x) = \frac{1}{(2\pi h)^d}\iint_{(\frac{x+y}{2},\theta) \in\Lambda_\Phi}e^{\frac{i}{h}(x-y)\cdot\theta}\frakp(\frac{x+y}{2},\theta) u(y)\,dy\wedge d\theta.
$$
Here
$$
\Lambda_\Phi = \left\{(x,\frac{2}{i}\partial_x \Phi(x))\::\:x\in\Bbb{C}^d\right\}
$$
with holomorphic gradient, and the natural symbol class for $\frakp:\Lambda_\Phi\rightarrow \Bbb{C}$ is
\begin{equation}
\label{eFBISymbol}
S(\Lambda_\Phi, \mathfrak{m}) := \{\mathfrak{a}\in C^{\infty}(\Lambda_\Phi,\Bbb{C})\::\:|\partial^\alpha_{x,\xi}\mathfrak{a}|\leq C_\alpha \mathfrak{m}\}
\end{equation}
for $\mathfrak{m}$ an appropriate order function.  When $\mathfrak{m}=1$, we obtain a uniformly bounded operator
$$
\operatorname{Op}^w_{\Phi,h}(\frakp) = \BigO(1): H_\Phi \rightarrow H_\Phi.
$$

The connection between the Weyl quantization on the real side and the same on the FBI-Bargmann transform side is made through the exact Egorov theorem,
\begin{equation}
\label{eEgorov}
T_\varphi p^w(x,hD_x) = \Op{\Phi}{h}(\frakp)T_\varphi,
\end{equation}
with
$$
\frakp \circ \varkappa_\varphi = p,
$$
\begin{equation}\label{defkap}
\varkappa_\varphi:(y, -\varphi'_y(x,y))\mapsto (x, \varphi'_x(x,y)).
\end{equation}
Using also (\ref{eqnPhiQuadraticDef}), it can be deduced that
$$
\Lambda_\Phi = \varkappa_\varphi(\Bbb{R}^{2d}).
$$
We note here that, with $\varphi_0$ in (\ref{ephi0}), we have
\begin{equation}\label{ekap0}
\kap_{\varphi_0}(y, \eta) = (y - i\eta, \eta).
\end{equation}

Where more rapid convergence of the integral is convenient, the contour of integration $\Lambda_\Phi$ may be replaced with
$$
\Gamma_{t_0} := \left\{\theta = \frac{2}{i}(\partial_x \Phi)(\frac{x+y}{2}) + it_0\overline{(x-y)}\right\},
$$
introducing an almost holomorphic extension of $\frakp$ to a tubular neighborhood of $\Lambda_{\Phi}$ and adding gaussian decay of the integrand off the main diagonal $\{x = y\}.$ In addition, a cutoff function of the form $\psi(x-y)$ for $\psi:\Bbb{C}^d\rightarrow[0,1]$ smooth and compactly supported with $\psi(x)=1$ near $x=0$, may be introduced into the integral. Both steps introduce an error term of the form
\begin{equation}\label{eqnFBICutoffContourShiftError}
R = \BigO(h^\infty): H_\Phi \rightarrow L^2_\Phi,
\end{equation}
and, modulo this error, we obtain the integral expression
\begin{equation}\label{eqnFBIQuantizationIntegralCutoff}
\Op{\Phi}{h}(\frakp)(u)(x) = \frac{1}{(2\pi h)^d} \iint_{\Gamma_{t_0}}e^{\frac{i}{h}(x-y)\cdot\theta}\psi_0(x-y)\frakp(\frac{x+y}{2},\theta) u(y)\,dy\wedge d\theta,
\end{equation}
where we continue to write $\frakp$ for an almost holomorphic extension.

Another view of the contour of integration $\Gamma_{t_0}$ is given by noting that that $|dy\wedge d\theta|$ pulls back
to a multiple of the Lebesgue volume form on $\Bbb{C}^d$,
\[
dL(y) = \bigwedge_{j=1}^d (d\jvRe y_j \wedge d\jvIm y_j) = 2^{-d} |dy\wedge d\bar{y}|.
\]
In fact, for quadratic $\Phi$ (thus with constant second derivatives), the definition of $\theta$ in $\Gamma_{t_0}$ above gives
\begin{eqnarray*}
|dy\wedge d\theta| &=& \left| dy\wedge\left(\frac{2}{i}(\partial_x^{2}\Phi)dy + (\frac{2}{i}\bar{\partial}_x\partial_x \Phi)d\bar{y} + (-it_0) d\bar{y}\right)\right|
\\ &=& \left|\operatorname{det}\left(2\bar{\partial}_x\partial_x \Phi + t_0 I\right)dy \wedge d\bar{y}\right|.
\end{eqnarray*}

A brief computation reveals that, if
\[
\varphi(x,y) = \frac{1}{2}\langle x,Ax \rangle + \langle x, By \rangle + \frac{1}{2}\langle y, Cy\rangle,
\]
then the critical value in the
definition of $\Phi(x)$ in (\ref{eqnPhiQuadraticDef}) is attained at
\[
y_0(x) = -(\jvIm C)^{-1} \jvIm(B^t x).
\]
Here $A$, $B$, and $C$ are complex matrices with $A$ and $C$ symmetric, $B$ invertible, and $\jvIm C$ positive definite. It can then be computed that
\[
\Phi(x) = -\frac{1}{2}\jvIm(\langle x, Ax \rangle) + \frac{1}{2}\langle \jvIm(B^t x), (\jvIm C)^{-1} \jvIm(B^t x)\rangle.
\]
To compute $\bar{\partial}_x \partial_x \Phi,$ note that, for holomorphic $f,$ we may write $\jvIm f$ as the sum of holomorphic and antiholomorphic parts, and, as a consequence, \[\bar{\partial}_x\partial_x \jvIm\langle x, Ax \rangle = 0.\]  Expanding $\jvIm(B^t
x) = \jvIm B^t \jvRe x + \jvRe B^t \jvIm x,$ multiplying out the second inner product term above, and using \[\bar{\partial}_x\partial_x = \frac{1}{4}(\partial_{\jvRe x}^2 +
\partial_{\jvIm x}^2),\] we see that
\[
\bar{\partial}_x\partial_x \Phi = \frac{1}{4}(\jvIm B (\jvIm C)^{-1} \jvIm B^t + \jvRe B(\jvIm C)^{-1} \jvRe B^t).
\]
Since $\operatorname{det} B \neq 0$ and $\jvIm C > 0,$ we recover the well-known fact that $\Phi$ is uniformly strictly plurisubharmonic, 
so that $\bar{\partial}_x\partial_x \Phi$ is a positive definite quadratic form. It follows that 
$|dy\wedge d\theta|$ is a constant non-zero multiple of $dL(y)$ for fixed $\varphi$ and $t_0 \geq 0$.

Finally, we note that, given $\varkappa$ a linear canonical map on $\Bbb{C}^{2d}$ sufficiently close to $\kap_{\varphi_0}$ in (\ref{ekap0}), we 
can easily obtain a $\varphi = \varphi(\varkappa)$, a holomorphic quadratic form of the type allowable in FBI-Bargmann transforms, for which
$$
\varkappa = \varkappa_\varphi.
$$
In fact, the differential equation which results for $\varphi$ is exact if and only if the linear transformation $\varkappa$ is canonical.  This is shown in the following lemma.

\begin{lemma}\label{lemPhiKappaExists}
A holomorphic linear map
\[
\kap(y,\eta) = \left(\begin{array}{cc}A & B \\ C & D\end{array}\right):\Bbb{C}^{2d}\rightarrow \Bbb{C}^{2d}
\]
is given by a holomorphic quadratic $\varphi_\kap:\Bbb{C}^{2d}\rightarrow \Bbb{C}$ as in (\ref{defkap}) if and only if $B$ is invertible and $\kap$ is canonical.  Both of these are equivalent to the three conditions,
\renewcommand{\labelenumi}{(\roman{enumi})}
\begin{itemize}
\item[(i)] $(DB^{-1})^t = DB^{-1},$
\item[(ii)] $(B^{-1}A)^t = B^{-1}A,$ and
\item[(iii)] $-(B^{-1})^t = C - DB^{-1}A.$
\end{itemize}

\end{lemma}

\begin{proof}

Begin by noting that any $\varphi:\Bbb{C}^{2d}\rightarrow\Bbb{C}$ holomorphic quadratic may be written as
\[
\varphi(x,y) = \frac{1}{2}\langle x, (\varphi_{xx}'') x\rangle + \langle x, (\varphi_{xy}'')y\rangle + \frac{1}{2}\langle y, (\varphi_{yy}'')y\rangle.
\]
Here, $(\varphi_{xy}'') = (\partial_{x_i}\partial_{y_j}\varphi)_{i,j=1}^n$ and inner products must be symmetric, meaning without complex conjugates.  Writing $(x,\xi) = \kap(y,\eta)$ for $\kap = \kap_\varphi$ as in (\ref{defkap}), we see that
\[
\varphi_y'(x,y) = (\varphi_{xy}'')^t x + (\varphi_{yy}'') y = -\eta(x,y).
\]
On the other hand,
\[
x = Ay + B\eta \Rightarrow \eta = B^{-1}x - B^{-1}Ay.\]  Therefore \[\varphi'_y = -B^{-1}x + B^{-1}Ay,
\]
from which we can deduce
\[
\varphi_{xy}'' = -(B^{-1})^t; \quad \varphi_{yy}'' = B^{-1}A.
\]
Similarly writing 
\[
\varphi_x' = Cy + D\eta(x,y) = DB^{-1}x + (C-DB^{-1}A)y,
\]
we get the final two relations
\[
\varphi_{xx}'' = DB^{-1};~~~ \varphi_{xy}'' = C-DB^{-1}A.
\]
Noting that the differential equation is exact if and only if $\varphi_{xx}''$ and $\varphi_{yy}''$ are symmetric and the formulas for $\varphi_{xy}''$ are equal gives the three conditions in the lemma.

Writing the complex symplectic inner product as
\[
\sigma\left((x,\xi), (y, \eta)\right) = \langle (x,\xi), J(y, \eta)\rangle = \left\langle \left(\begin{array}{c} x\\ \xi \end{array} \right), \left(\begin{array}{cc} 0 & -1 \\ 1 & 0 \end{array}\right)\left(\begin{array}{c} y \\ \eta \end{array}\right)\right\rangle,
\]
we see that $\kap$ is canonical if and only if $\kap^t J \kap = J,$ which by computation is equivalent to
\[
\left(\begin{array}{cc} -A^t C + C^t A & -A^t D + C^t B \\ -B^t C + D^t A & -B^t D + D^t B\end{array}\right) = \left(\begin{array}{cc} 0 & -1 \\ 1 & 0\end{array}\right).
\]
This encodes three equalities, given by the upper-left, lower-right, and upper-right entries of the block matrix.  (Upper-right and lower-left are transparently equivalent.)

Upon assuming that $B$ is invertible, the lower-right equality is equivalent to (i) in the lemma by premultiplying by $(B^{-1})^t$ and postmulitplying by $B^{-1}.$  When (i) holds, we have that the upper-right equality and (iii) in the theorem are equivalent after premultiplying by $(B^{-1})^t.$  Finally, using (iii) in the lemma to replace $C$ in the upper-left equality and cancelling $A^tDB^{-1}A$ gives equivalence with (ii), concluding the proof that the conditions are equivalent.

\end{proof}

\subsection{Shifting to weights near $\Phi_0$}\label{ssNearWeights}

Associated to the weight function $G_{\eps}$, whose properties were reviewed in proposition \ref{propGeps}, we introduce the IR-manifold
\begin{equation}
\label{eLambdaDeltaReal}
\Lambda_{\delta,\eps} = \left\{X+i\delta H_{G_{\eps}}(X)\::\:X \in \Bbb{R}^{2d}\right\},
\end{equation}
defined for $\delta>0$ and $\eps>0$ small enough. In what follows, the small number $\delta>0$ will be kept fixed, and the dependence on $\delta$
in estimates will therefore not be indicated explicitly.

Arguing as in \citep{HeSjSt05} and \citep{H-PS}, we obtain that if we define
\begin{equation}
\label{eqnDefPhiEps}
\Phi_\eps(x) = \operatorname{v.c.}_{(y, \eta)\in \Bbb{C}^d \times \Bbb{R}^d}\left(-\jvIm \varphi_0(x,y) -
(\jvIm y)\cdot \eta + \delta G_\eps(\jvRe y, \eta)\right),
\end{equation}
where $\operatorname{v.c.}$ stands for the critical value, then
\begin{equation}
\label{eLambdaDeltaFBI}
\kappa_{\varphi_0}\left(\Lambda_{\delta,\eps}\right)=
\Lambda_{\Phi_\eps} := \left\{\left(x, \frac{2}{i}\partial_x \Phi_\eps (x) \right)\::\: x \in \Bbb{C}^d\right\}.\end{equation}
Continuing to follow \citep{HeSjSt05}, section 3, one can check that $\Phi_{\eps}\in C^{\infty}(\Bbb{C}^{d})$ is a uniformly strictly plurisubharmonic function such that
\[
\Phi_\eps(x) = \Phi_0(x) + \delta G_{\eps}(\jvRe x, -\jvIm x) + \BigO(\delta^2\eps).
\]
Furthermore, $\Phi_{\eps}-\Phi_0$ is compactly supported and we have the following basic properties, valid uniformly in $\eps>0$:
\begin{equation} \label{eqnPhisClose}
||\Phi_\eps-\Phi_0||_{L^\infty} = \BigO(\eps),
\end{equation}
\begin{equation}\label{eqnPhiDerivsClose}
||\nabla(\Phi_\eps - \Phi_0)||_{L^\infty} = \BigO(\eps^{1/2}),
\end{equation} and
\begin{equation} \label{eqnPhiHessianBounded}
||\nabla^2\Phi_\eps||_{L^\infty} = \BigO(1).
\end{equation}

For future reference, let us now recall the linear IR-manifold $\Lambda_{\delta}$, introduced in (\ref{eqnDefLambdaDelta}). Following \citep{H-PS}, we then find that
$$
\kap_{\varphi_0}(\Lambda_{\delta}) = \Lambda_{\Phi_q} = \left\{\left(x, \frac{2}{i}\partial_x \Phi_q(x) \right)\::\: x \in \Bbb{C}^d\right\},
$$
where $\Phi_q$ is a strictly plurisubharmonic quadratic form on $\Bbb{C}^d$, satisfying
\begin{equation}\label{ePhiqClose}
\Phi_q(x) = \Phi_0(x) + \delta G_q(\jvRe x, -\jvIm x) +\BigO(\delta^2 |x|^2).
\end{equation}
The quadratic weight function $\Phi_q$ can also be given as a critical value, similarly to (\ref{eqnDefPhiEps}), or via 
lemma \ref{lemPhiKappaExists} applied to $\kap_{\varphi_0}\circ\kap_q$ of (\ref{ekapq}). We shall also have to recall,
following \citep{H-PS}, that in a tiny neighborhood of zero, the weight $\Phi_{\eps}$ is close to $\Phi_q$, and, specifically,
\begin{equation}
\label{q-eps}
\Phi_{\eps}(x) = \Phi_q(x) + \BigO(|x|^3),~~~\forall |x| \leq \eps^{1/2}.
\end{equation}
The estimate (\ref{q-eps}) will be important in making estimates localized to a neighborhood of size $|x| \leq \eps^{1/2}.$
Let us also remark that we will henceforth consider the time $T$ from (\ref{eqnDefGq}), which is also implicitly in proposition \ref{propGeps},
as well as $\delta > 0$, to be fixed.
While there is some restriction on these constants, we may choose $T$ and $\delta$ sufficiently small that, after considerations in proposition \ref{propGeps} and the asymptotic expansion for $\Phi_\eps,$ we may leave them $h$-independent.

Associated to the function $\Phi_{\eps}$ is the weighted space $H_{\Phi_\eps, h}$, and the $L^\infty$ bound (\ref{eqnPhisClose}) gives upper and lower bounds for the norm of the ``identity map"
\[
H_{\Phi_0, h} \ni u \mapsto u \in H_{\Phi_\eps, h}.
\]
The weights attached to $dL(x)$ are then governed by the ratio
\[
e^{-C\frac{\eps}{h}} \leq \frac{e^{-\frac{2}{h}\Phi_\eps(x)}}{e^{-\frac{2}{h}\Phi_0(x)}} \leq e^{C\frac{\eps}{h}},
\]
for some $C > 0,$ and so the norms obey
\begin{equation}\label{eqnNormsClose}
e^{-C\frac{\eps}{h}} \leq \frac{||u||_{H_{\Phi_\eps}}}{||u||_{H_{\Phi_0}}} \leq e^{C\frac{\eps}{h}}.
\end{equation}
Note that this also implies that, as subsets of $\operatorname{Hol}(\Bbb{C}^d),$ the normed spaces $H_{\Phi_\eps}$ and $H_{\Phi_0}$ are identical, with the norms equivalent for each fixed $h>0$ but not uniformly as $h\rightarrow 0$. In order to shift between the various weighted spaces with losses limited to a negative power of $h$, we will use
\begin{equation}
\label{eqnEpsVagueDef}
\varepsilon = \frac{1}{C}h\log\frac{1}{h},
\end{equation}
for $C > 0$ to be chosen.  As a consequence, the ratio in (\ref{eqnNormsClose}) is bounded by 
\[
	e^{\pm \BigO(\eps)/h} = h^c, c > 0.
\]
Thus, taking $C$ large will allow us to lose arbitrarily small fractional powers of $h.$

Working in the $H_{\Phi_{\eps}}$-spaces, with $\eps$ satisfying (\ref{eqnEpsVagueDef}), allows us to preserve much of the framework of 
the FBI transform introduced in section \ref{subsecFBIQuadraticPhi}. As in (\ref{eqnFBICutoffContourShiftError}), we have errors which are $\BigO(h^{\infty}):L^2_{\Phi_0}\rightarrow L^2_{\Phi_0},$ and those errors with respect to $\Phi_\eps$ become
\begin{equation} \label{error}
\BigO(e^{C\eps/h}h^\infty) =
\BigO(h^\infty):L^2_{\Phi_\eps}\rightarrow L^2_{\Phi_\eps},
\end{equation}
for $\eps$ as in (\ref{eqnEpsVagueDef}).

To study the small error introduced by replacing $\Op{\Phi_0}{h}(\frakp)$ with $\Op{\Phi_\eps}{h}(\frakp)$ for $\frakp$ extended 
almost analytically off of $\Lambda_{\Phi_0}$, we introduce a parameterized family of contours via
\begin{equation} \label{eqnGammat}
\Gamma_{t,t_0} := \left\{\theta = \frac{2}{i}\partial_x\left((1-t)\Phi_0 + t\Phi_\eps\right)(\frac{x+y}{2}) + it_0\overline{(x-y)}\right\},
\end{equation}
for $t \in [0,1]$, and $t_0$ large and fixed.  On the real side, the fact that the extension of $p$ off $\Bbb{R}^{2d}$ is almost analytic means that
\[
\bar{\partial}_{z,\zeta}p(z,\zeta) = \BigO(|\jvIm(z,\zeta)|^\infty).
\]
As $\kap_{\varphi_0}$ is linear, it is not hard to show that this implies that $\frakp = p \circ \kap_{\varphi_0}^{-1}$ is almost analytic in the analogous sense that
\begin{equation}\label{eFBIAlmostHol}
\bar{\partial}_{x,\xi}\frakp(x,\xi) = \BigO\left(\operatorname{dist}\left((x,\xi),\Lambda_{\Phi_0}\right)^{\infty}\right).
\end{equation}

Then, using Stokes' formula, we obtain for $u\in H_{\Phi_{\eps}}$, neglecting error terms of the form (\ref{error}),
\begin{multline*}
\Op{\Phi_\eps}{h}(\frakp)u(x) = \frac{1}{(2\pi h)^d}\iint_{\Gamma_{1,t_0}}e^{\frac{i}{h}(x-y)\cdot\theta}\psi_0(x-y)\frakp(\frac{x+y}{2},\theta) u(y)\,dy\wedge d\theta
\\ + (2\pi h)^{-d} \iiint_{\Gamma_{[0,1]}} e^{\frac{i}{h}(x-y)\cdot\theta} \bar {\partial}_{y,\theta}\left(\psi_0(x-y) \frakp(\frac{x+y}{2},\theta)\right) \wedge dy \wedge d\theta,
\end{multline*}
where $\Gamma_{[0,1]}$ is the natural union of the
$\Gamma_{t,t_0}$ of (\ref{eqnGammat}), where $0\leq t\leq 1$.  The first integral defines a uniformly bounded operator on $L^2_{\Phi_\eps},$ and
details of the estimate on the phase may be found in the proof of proposition \ref{pLocalSmallErrors} below. Along $\Gamma_{[0,1]}$, by
(\ref{eqnPhiDerivsClose}) and (\ref{eFBIAlmostHol}),
\[
\bar{\partial}_{y,\theta} \left(\psi_0(x-y)\frakp(\frac{x+y}{2},\theta)\right) = \BigO(|x-y| + \eps^{1/2})^\infty.
\]
In order to estimate the second term as a map on $H_{\Phi_{\eps}}$, we use Schur's test together with the uniform estimates (\ref{eqnPhisClose}) and (\ref{eqnPhiHessianBounded}), which allows us to bound the effective kernel of the second term by an expression of the form
$$
e^{\frac{C\eps}{h}} h^{-d} e^{-t_0|x-y|^2/h}\BigO\left((|x-y| + \eps^{1/2})^\infty\right).
$$
The corresponding contribution to $\Op{\Phi_\eps}{h}({\frakp})$ is therefore
$$
e^{C\eps/h}\BigO(h^\infty + \eps^\infty),
$$
and since $\eps$ is only logarithmically larger than $h$, we conclude that, when $u\in H_{\Phi_{\eps}}$,
\[
\Op{\Phi_\eps}{h}(\frakp)= \frac{1}{(2\pi
h)^d}\iint_{\Gamma_{1,t_0}}e^{\frac{i}{h}(x-y)\cdot\theta}\psi_0(x-y)\frakp(\frac{x+y}{2},\theta) u(y)\,dy\wedge
d\theta + Ru.
\]
Here the operator norm of $R$, viewed as a map from $L^2_{\Phi_{\eps}}$ to itself, is $\BigO(h^\infty)$.

\medskip
\noindent\emph{Remark}.  When comparing spaces of holomorphic functions which are $L^2$ against differing exponential weights, 
it is natural to ask whether the spaces, as linear subspaces of the set of holomorphic functions, share any elements at all.  We consider our current situation, where $\Phi_0$ may be compared with a quadratic $\Phi_q$ obeying (\ref{ePhiqClose}) for $\delta > 0$ small.

A convenient place to look is the images under $T_{\varphi_0}$ of eigenfunctions $\{v_j\}$ of an elliptic quadratic form, e.g. the harmonic oscillator
\[
q_0(x,hD_x) = x^2 + (hD_x)^2.
\]
This is because the eigenfunctions form a basis of $L^2(\Bbb{R}^d)$ and, by the exact Egorov theorem, we should expect
$T_{\varphi_\delta}^*T_{\varphi_0}v_j$ to be eigenfunctions of $q_0 \circ \kappa_{\varphi_0}\circ \kappa_{\varphi_\delta},$ which remains
elliptic when $\delta$ is small.

Indeed, the eigenfunctions of $q_0(x,hD_x)$ are given by an algebraic basis of $P(\Bbb{R}^d)e^{-x^2/2h},$ and a simple computation shows
that
\[
T_{\varphi_0}(P(\Bbb{R}^d)e^{-y^2/2h}) = P(\Bbb{C}^d)e^{-x^2/4h}.
\]
For a (holomorphic) polynomial $f(x),$ we have
\[
f(x)e^{-x^2/4h} = \BigO_f(1)e^{-c|x|^2/h}e^{\Phi_0(x)/h}
\]
for any $c < 1/4.$  This additional convergence factor allows us to conclude by Stokes' theorem that, for $u \in P(\Bbb{C}^d)e^{-x^2/4h},$
\[
\Op{\Phi_q}{h}(\frakq)(u) = \Op{\Phi_0}{h}(\frakq)u
\]
and we see that $P(\Bbb{C}^d)e^{-x^2/4h}$ is dense in both $H_{\Phi_0}$ and $H_{\Phi_q}$ when $\delta > 0$ is sufficiently small.

We may also recall, following \citep{Sj74}, that the generalized eigenfunctions of any quadratic differential operator $q^w(x,D_x)$ such that
$\jvRe q(x,\xi)\geq |(x,\xi)|^2/C$ are given by
\begin{equation}
\label{eqEig}
p(x) e^{\Phi(x)},
\end{equation}
where $p(x)$ is a polynomial and $\Phi(x)$ is a quadratic form with $\jvIm \Phi(x)\geq |x|^2/C$.  Let us also recall that $\Phi$ in (\ref{eqEig})
is such that the positive Lagrangian subspace $\{(x,\partial_x \Phi(x));\, x\in \Bbb{C}^d\}$ is the direct sum of the generalized eigenspaces of the Hamilton map of $q$, corresponding to the eigenvalues with positive imaginary part.

\subsection{Rescaling}\label{subsecRescaling}

We use FBI-side changes of variables
\begin{equation}\label{FBIChgVarsDef}
\frakU_\varepsilon: H_{\Phi}(\Bbb{C}^d;h)\ni u(x) \mapsto \varepsilon^{d/2}u(\varepsilon^{1/2}x)\in H_{\tilde{\Phi}}(\Bbb{C}^d; \tilde{h}).
\end{equation}
With
$$
\tilde{h} := \frac{h}{\varepsilon}, \quad \tilde{\Phi}(x) := \frac{1}{\varepsilon}\Phi(\varepsilon^{1/2}x),
$$
this change of variables is unitary. Furthermore, writing 
\[
	\frakp_\varepsilon(x,\xi) = \frakp(\varepsilon^{1/2}x,\varepsilon^{1/2}\xi),
\]
the change of variables interacts with quantizations via the property
\begin{equation}
\label{FBIChgVarsSymbol}
\frakU_\varepsilon \Op{\Phi}{h}(\frakp) = \Op{\tilde{\Phi}}{\tilde{h}}(\frakp_\varepsilon)\frakU_\varepsilon.
\end{equation}

The natural real-side analogues are the operator
\begin{equation}
\label{eRealRescaling}
U_\varepsilon:L^2(\Bbb{R}^d)\ni v(y)\mapsto \varepsilon^{d/4}v(\varepsilon^{1/2}y)\in L^2(\Bbb{R}^d),
\end{equation}
which is automatically unitary, and the symbol transformation rule
$$
U_\varepsilon p^w(y,hD_y) = p_\varepsilon^w(y, \tilde{h}D_y) U_\varepsilon
$$
with $p_\varepsilon(x,\xi) = p(\eps^{1/2}x, \eps^{1/2}\xi).$

In cases where the symbol or the weight is quadratic, the rescaling is much simpler.  On the real side with $q(x,\xi)$ a quadratic form, we have
\begin{equation}
\label{RealQuadraticScaling}
U_\varepsilon q^w(y, hD_y) = \varepsilon q^w(y, \tilde{h}D_y)U_\varepsilon,
\end{equation}
and if $\Phi:\Bbb{C}^d\rightarrow \Bbb{C}$ is quadratic (as well as $\frakq:\Bbb{C}^{2d}\rightarrow \Bbb{C}$), then $\tilde{\Phi} = \Phi$ and so
\begin{equation}
\label{FBIQuadraticScaling}
\frakU_\varepsilon \Op{\Phi}{h}(\frakq) = \varepsilon\Op{\Phi}{\tilde{h}}(\frakq)\frakU_\varepsilon.
\end{equation}

\subsection{Quantization-multiplication formula}

To be able to handle regions of the phase space on which the symbol of the operator is sufficiently elliptic, we shall
make use of a basic formula relating the action of the operator and multiplication by the symbol, on the level of inner products.

\begin{proposition} Let $\psi(x) \in C_b^\infty(\Bbb{C}^d;[0,1])$ be such that $\nabla\psi$ is compactly supported. Assume that 
$\frakp\in C^{\infty}(\Lambda_{\Phi_{\eps}})$ is an almost holomorphic extension of a symbol on $\Lambda_{\Phi_0}$ obeying, uniformly in $\eps>0$ along $\Lambda_{\Phi_{\eps}}$,
$$
|\partial_{x,\xi}^\alpha p(x,\xi)| = \BigO_\alpha(1), \quad \forall |\alpha| \geq 2.
$$
Writing $\xi(x) = \frac{2}{i} \partial_x \Phi_\eps(x)$, the quantization-multiplication formula
\begin{equation}\label{QM}
\langle \psi \Op{{\Phi_{\eps}}}{h}(\frakp)u, u\rangle_{H_{\Phi_{\eps}}} =
\int_{\Bbb{C}^d} \psi(x) \frakp(x,\xi(x))|u(x)|^2 e^{-\frac{2}{h}\Phi_{\eps}(x)}\,dL(x) + \BigO(h)||u||^2_{H_{\Phi_{\eps}}}
\end{equation}
holds for all $u\in H_{\Phi_{\eps}}$ with error a function in $L^2_{\Phi_\eps}$.
\end{proposition}

Various forms of this formula are proven in \citep{HeSjSt05}, \citep{H-PS}, \citep{SjLoR}, \citep{Sj90}, and \citep{Vi09}, among others.

\begin{proof} Begin by Taylor expanding $\frakp$ along
$\Lambda_{\Phi_{\eps}}$ at the point $(x,\xi(x))$ to obtain approximate
values for $\frakp$ at $((x+y)/2, \theta):$
\begin{multline}
\label{QMTaylor}
\frakp(\frac{x+y}{2},\theta) = \frakp(x,\xi(x)) +
\sum_{j=1}^d (\partial_{\theta_j}\frakp)(x,\xi(x))(\theta_j - \xi_j(x)) \\
+ \sum_{j=1}^d (\partial_{x_j}\frakp)(x,\xi(x))(\frac{y_j-x_j}{2}) +
\mathfrak{r}(x,y,\theta),
\end{multline} with
\[
|\mathfrak{r}(x,y,\theta)|\leq ||\nabla^2\frakp||_{L^\infty}\left(|\theta - \xi(x)|^2 +
\frac{|y-x|^2}{4}\right)+\mathcal{O}(h^{\infty}).
\]
Here we have used holomorphic derivatives, relying on the fact that we are working with an almost holomorphic extension off $\Lambda_{\Phi_0}$, and the $\BigO(h^\infty)$ error in $\mathfrak{r}$ comes from the distance between $\Lambda_{\Phi_0}$ and $\Lambda_{\Phi_\eps}$.  We have seen in section \ref{ssNearWeights} that 
$\Op{{\Phi_{\eps}}}{h}(\frakp)$ is realized along the contour
$$
\Gamma_{t_0}: \theta  = \frac{2}{i} \frac{\partial \Phi_{\eps}}{\partial x} \left(\frac{x+y}{2}\right)+it_0\overline{(x-y)},\quad t_0>0,
$$
and therefore, we conclude from the definition of the contour that
\[
|\theta - \xi(x)| \leq \BigO(||\nabla^2\Phi||_{L^\infty})|x-y|,
\]
so that, uniformly in $\eps>0$,
\[
\mathfrak{r} = \BigO(|x-y|^2 + h^{\infty})
\]
on the contour in the definition of $\Op{{\Phi_{\eps}}}{h}(\frakp).$ Neglecting $\BigO(h^{\infty})$ errors, it follows that the effective integral kernel for
$\Op{{\Phi_{\eps}}}{h}(\mathfrak{r})$ on $L^2_{\Phi_{\eps}}$ is bounded by
\[
\BigO(h^{-d}|x-y|^2 e^{-\frac{C}{h}|x-y|^2}) = h \BigO(h^{-d}\frac{|x-y|^2}{h}e^{-\frac{C}{h}|x-y|^2}).
\]
Schur's test then implies that
\[
\Op{{\Phi_{\eps}}}{h}(\mathfrak{r}) = \BigO(h):L^2(\Bbb{C}^d, e^{-2\Phi_{\eps}/h}\,dL(x))\rightarrow
L^2(\Bbb{C}^d, e^{-2\Phi_{\eps}/h}\,dL(x)),
\]
and so, by the Cauchy-Schwarz inequality,
\[
	\langle{\psi \Op{{\Phi_{\eps}}}{h}(\mathfrak{r})u, u\rangle}_{H_{\Phi_{\eps}, h}} = \BigO(h)||u||^2_{H_{\Phi_{\eps}, h}}.
\]

The closed contour formed by the difference between
$\Lambda_{\Phi_0}$ and $\Lambda_{\Phi_{\eps}}$ is bounded as a consequence of
$\Phi_0 = \Phi_{\eps}$ for sufficiently large $|x|.$  Therefore, a bounded change of contour with holomorphic integrand establishes the formulas
\[
(2\pi h)^{-d}\iint_{(\frac{x+y}{2}, \theta)\in \Lambda_{\Phi_{\eps}}} y_je^{\frac{i}{h}(x-y)\cdot\theta}\,dy\,d\theta = x_ju(x)
\]
and
\[
(2\pi h)^{-d}\iint_{(\frac{x+y}{2}, \theta)\in \Lambda_{\Phi_{\eps}}} \theta_j e^{\frac{i}{h}(x-y)\cdot\theta}\,dy\,d\theta = hD_{x_j}u(x),
\]
as a consequence of the standard formulas, where the contour is along $\Lambda_{\Phi_0}.$

These formulas allow us to simplify the integral for the second and third parts of the Taylor expansion (\ref{QMTaylor}):
\begin{multline*}
(2\pi h)^{-d}\sum_{j=1}^d(\partial_{x_j}\frakp)(x,\xi(x))\iint_{\Gamma_{t_0}} \frac{y_j-x_j}{2}e^{\frac{i}{h}(x-y)\cdot\theta}u(y)\,dy\,d\theta
\\ = \sum_{j=1}^d (\partial_{x_j}\frakp)(x,\xi(x))\frac{x_j - x_j}{2}u(x) = 0,\end{multline*} and \begin{multline*}(2\pi h)^{-d}\sum_{j=1}^d (\partial_{x_j}\frakp)(x,\xi(x))\iint_{\Gamma_{t_0}} (\theta_j - \xi_j(x))e^{\frac{i}{h}(x-y)\cdot\theta}u(y)\,dy\,d\theta
\\ = \sum_{j=1}^d (\partial_{\theta_j}\frakp)(x,\xi(x))(hD_{x_j} - \xi_j(x))u(x).\end{multline*}

The contribution from $(\partial_{\theta_j}\frakp)(x,\xi(x))(hD_{x_j} - \xi_j(x))u(x)$ may be bounded by integration by parts
$\int (D_{x_j}f)g = -\int f(D_{x_j}g)$ applied with $f = u(x)$ and
$g = \psi(x)(\partial_{\theta_j}\frakp)(x,\xi(x))\overline{u(x)}e^{-2\Phi_{\eps}(x)/h}.$
By the definition of $\xi(x)$,
$$
(-hD_{x_j} - \xi_j) e^{-\frac{2}{h}\Phi(x)} = 0.
$$
Thus canceling the term where $hD_x$ hits the exponential weight, we are left with
\begin{multline*}
\langle \psi(x) (\partial_{\theta_j}\frakp)(x,\xi(x)) (hD_{x_j} - \xi_j(x))u(x), u(x)\rangle_{H_{\Phi_{\eps}}}
\\ = -\int_{\Bbb{C}^d}u(x) (hD_{x_j})\left(\psi(x)(\partial_{\theta_j}\frakp) (x,\xi(x))\overline{u(x)}\right) e^{-\frac{2}{h}\Phi_{\eps}(x)} \,dL(x).
\end{multline*}
The function $\overline{u(x)}$ is antiholomorphic and therefore commutes with $D_{x_j}.$  The derivative hitting $\frakp$ is controlled by
\[
|D_x(\partial_{\theta_j}\frakp(x,\xi(x)))| \leq ||\nabla^2\frakp||_{L^\infty}||\nabla \xi||_{L^\infty} = ||\nabla^2\frakp||_{L^\infty}||\nabla^2\Phi||_{L^\infty} = \BigO(1).
\]
Finally, since $\nabla \psi$ is compactly supported, the contribution of the term containing this function is clearly harmless. The proof is complete.

\end{proof}

\section{Improved resolvent estimates for quadratic operators}\label{sImprovedQuad}

In the author's previous work \citep{Vi09}, Proposition 3.1, trace-class perturbations for general pseudodifferential operators were used to obtain resolvent estimates, of the form
$$
	||(q^w(x,hD_x) - z)^{-1}||_{L^2(\Bbb{R}^d)\rightarrow L^2(\Bbb{R}^d)} \leq \BigO(h^{-1-\gamma}),
$$
in the case where $q:\Bbb{R}^{2d}\rightarrow \Bbb{C}$ is an elliptic quadratic form.  The spectral parameter $z \in \Bbb{C}$ was restricted to the region $|z|\leq hf(h)$ and assumed to obey
$$
	\opnm{dist}(z, \opnm{Spec}(q^w(x,hD_x))) \geq hf(h)^{(1-d)/2},
$$
for
$$
	f(h) = \frac{1}{M}\left(\frac{\log \frac{1}{h}}{\log\log\frac{1}{h}}\right)^{1/d},
$$
with the constant $M$ large but fixed based on $\gamma > 0$ and $q$.

Using a more careful analysis of the lattice structure (\ref{eqnSpectrumLattice}) of the eigenvalues for $q^w(x,hD_x)$, one may obtain the following improvement, both in $f(h)$ and in rapid approach to the spectrum, of this estimate.

\begin{proposition}\label{pImprovedQuadratic}

Fix $\gamma > 0.$  Let $q$ be a quadratic form on $\Bbb{R}^{2d}$, elliptic in the sense that 
\begin{equation}\label{eIQell}
	\jvRe q(x,\xi) \geq \frac{1}{C} |(x,\xi)|^2.
\end{equation}
Define
\begin{equation}\label{eIQf}
	f(h) = \frac{1}{M}\left(\log\frac{1}{h}\right)^{1/d},
\end{equation}
for $M$ sufficiently large depending on $\gamma$ and $q.$  Then there exist some $C_0 > 0$ sufficiently large and some $h_0 > 0$ sufficiently small where, for any $h \in (0, h_0]$ and for any $z \in \Bbb{C}$ with $|z| \leq hf(h)$ and
\begin{equation}\label{eIQSpecDist}
	\opnm{dist}(z, \opnm{Spec}(q^w(x,hD_x))) \geq he^{-f(h)/C_0},
\end{equation}
we have the resolvent estimate
$$
	||(q^w(x,hD_x) - z)^{-1}||_{L^2(\Bbb{R}^d)\rightarrow L^2(\Bbb{R}^d)} \leq \BigO(h^{-1-\gamma}).
$$

\end{proposition}

\begin{proof}

Let us begin by introducing the notation
$$
	S_h(R) = \opnm{Spec}(q^w(x,hD_x))\cap \{|z|\leq R\},
$$
emphasizing that the semiclassical parameter $h$ may change.  From (\ref{RealQuadraticScaling}) it is easy to see that $S_h(R) = h S_1(h^{-1}R)$.

From (\ref{eIQell}) we see that $q(x,\xi)$ is elliptic near infinity in the symbol class $S(m)$, where 
$$
	m(x,\xi) = 1+\frac{1}{4C}|(x,\xi)|^2,
$$
and because $\jvRe q(x,\xi) \geq 0$ for all $(x,\xi)$. Therefore we may apply Proposition 2.1 from \citep{Vi09} with $\rho = 1$ and $\rho' = 2$, obtaining the semiclassical resolvent bound
$$
	||(q^w(x,\tilde{h}D_x) - \tilde{z})^{-1}||_{L^2(\Bbb{R}^d)\rightarrow L^2(\Bbb{R}^d)} \leq e^{C_1 \tilde{h}^{-d}}\prod_{\tilde{z}_j\in S_{\tilde{h}}(2)} |\tilde{z}-\tilde{z}_j|^{-1}.
$$
This holds for all $\tilde{z}$ with $|\tilde{z}|\leq 1$ in the limit $\tilde{h}\rightarrow 0^+.$  We will use
$$
	\tilde{h} = \frac{1}{f(h)}.
$$
Furthermore, write
$$
	\tilde{z} = \frac{z}{hf(h)}
$$
and recall that, from the change of variables (\ref{RealQuadraticScaling}), we have the unitary equivalence 
\[
	q^w(x,\tilde{h}D_x) \sim \frac{\tilde{h}}{h}q^w(x,hD_x).
\]
This provides the rescaled estimate
$$
	||(q^w(x,hD_x)-z)^{-1}|| \leq \frac{1}{hf(h)} e^{C_1\tilde{h}^{-d}}\prod_{\tilde{z}_j\in S_{\tilde{h}}(2)} |\tilde{z}-\tilde{z}_j|^{-1}.
$$
With $f(h)$ chosen as in (\ref{eIQf}), clearly
\begin{equation}\label{eIQfSmall}
	f(h) \ll h^{-\gamma/3}.
\end{equation}
It is therefore sufficient to show that
\begin{equation}\label{eIQexp}
	e^{C_1\tilde{h}^{-d}} \leq \BigO(h^{-\gamma/3})
\end{equation}
and
\begin{equation}\label{eIQprod}
	\prod_{\tilde{z}_j\in S_{\tilde{h}}(2)} |\tilde{z}-\tilde{z}_j|^{-1} \leq \BigO(h^{-\gamma/3}).
\end{equation}

The requirement (\ref{eIQexp}) necessitates the choice of $f(h)$.  Since 
\[
	\jvexp(C_1 \tilde{h}^{-d}) = \jvexp(C_1f(h)^d),
\]
taking logarithms shows that $f(h)$ defined by (\ref{eIQf}) with $M = 2C_1/\gamma$ is sufficient and necessary to establish (\ref{eIQexp}).  We next consider the spectrum in the product appearing in (\ref{eIQprod}).

Let $F = \frac{1}{2}H_q$ be the fundamental matrix of $q$.  We write $\lambda_1,\dots,\lambda_n$ for the eigenvalues of $F$ with $\jvIm\lambda_j > 0$, counted for algebraic multiplicity, and we recall from Section 3 of \citep{Sj74} that there are $n$ such.  Furthermore, write $\mu_j = \lambda_j/i$.  We introduce the notation, for $\Bff{x}\in\Bbb{R}^d$,
\begin{equation}\label{eIQmu}
	\mu(\Bff{x}) = \sum_{j=1}^n (1+2x_j)\mu_j.
\end{equation}
With this we obtain the convenient formula
\begin{equation}\label{eMuSpectrum}
	\operatorname{Spec}(q^w(x,D_x)) = \left\{\mu(\Bff{k})\::\: \Bff{k} \in (\Bbb{N}\cup\{0\})^d\right\}
\end{equation}
(cf. (\ref{eqnSpectrumLattice}), also from \citep{Sj74}).  As we must consider multiplicity in the spectrum, we will regard the set on the right with multiplicity as well.  The two multiplicites agree in that the algebraic multiplicity of $\lambda \in \operatorname{Spec}(q^w(x,D_x))$ is equal to the number of $\Bff{k}$ with $\mu(\Bff{k}) = \lambda$.

We remove the semiclassical dependence on the parameter $\tilde{h}$ by using the change of variables (\ref{RealQuadraticScaling}) on the left-hand side of (\ref{eIQprod}) to turn $S_{\tilde{h}}(2)$ into $S_1(2f(h))$:
\begin{equation}\label{eIQzetaProd}
	\prod_{\tilde{z}_j\in S_{\tilde{h}}(2)} |\tilde{z}-\tilde{z}_j|^{-1} = \prod_{\zeta_j \in S_1(2f(h))}  \tilde{h}^{-1}|\frac{z}{h} - \zeta_j|^{-1}.
\end{equation}
We will divide the $\zeta_j\in S_1(2f(h))$ into strips parallel to the imaginary axis of size $\sim 1$, and we begin by counting the $\zeta_j$ in such a strip.

We will now see that
\begin{equation}\label{eIQLatticeCount}
	\#\{\Bff{k}\in(\Bbb{N}\cup\{0\})^d\::\:|\rho - \jvRe\mu(\Bff{k})|\leq r\} = \BigO_{r,q}(f(h)^{d-1}),
\end{equation}
uniformly for $|\rho| \leq 3f(h)$, so long as $f(h)$ is sufficiently large, or equivalently, so long as $\tilde{h} \in (0, \tilde{h}_0]$ for $\tilde{h}_0 > 0$ sufficiently small.  The hypothesis $|\rho| \leq 3f(h)$ may be replaced by $|\rho| \leq Cf(h)$ for any fixed $C$, but $C = 3$ suffices here.

This is a straightforward consequence of the volume of a $d$-dimensional simplex.  To aid in the exposition, for $\Bff{k} = (k_1,\dots,k_d) \in \Bbb{Z}^d$, we define the box in $\Bbb{R}^d$ with corner at $\Bff{k}$ via the formula
$$
	B(\Bff{k}) = \{\Bff{x} = (x_1,\dots,x_d)\in\Bbb{R}^d\::\: x_j \in (k_j, k_j+1], j = 1,\dots, d\}.
$$
This is so that, for any $K \subseteq \Bbb{Z}^d$, we have
\begin{equation}\label{eIQKB}
	\# K = \opnm{vol}\left(\bigcup_{\Bff{k}\in K} B(\Bff{k})\right).
\end{equation}
Say that $\Bff{x} \in B(\Bff{k})$ for some $\Bff{k} \in (\Bbb{N}\cup\{0\})^d$ obeying
$$
	|\rho - \jvRe \mu(\Bff{k})| \leq r.
$$
Clearly, $\Bff{x} \in \Bbb{R}^d_+$ as $k_j \geq 0$ for all $j$.  Furthermore, from the definition (\ref{eIQmu}) of $\mu(\Bff{x})$ and the fact that $\jvRe \mu_j > 0$, we see that
$$
	\rho-r < \jvRe \mu(\Bff{x}) \leq \rho+r+\sum_{j=1}^d \jvRe \mu_j.
$$
Using the definition of $\mu(\Bff{x})$ once more, we see that
$$
	\rho - r - \sum_{j=1}^d \jvRe \mu_j < \sum_{j=1}^d 2x_j\jvRe \mu_j \leq \rho+r.
$$

Writing
$$
	T(R) = \{\Bff{x}\in \Bbb{R}^d_+ \::\: \sum_{j=1}^d 2x_j \jvRe \mu_j \leq R\},
$$
it is now clear from (\ref{eIQKB}) that to prove (\ref{eIQLatticeCount}) it suffices to bound
$$
	\opnm{vol}(T(\rho+r)) - \opnm{vol}(T(\rho-r-\sum_{j=1}^d \jvRe \mu_j)).
$$
Elementary change of variables and recalling that $\jvRe \mu_j > 0$ for all $j$ gives that
$$
	\opnm{vol}(T(R)) = \frac{1}{2^d d!}\left(\prod_{j=1}^d (\jvRe \mu_j)^{-1}\right)R^d = C_q R^d.
$$
This allows us to conclude that
$$
	\opnm{vol}(T(\rho+r)) - \opnm{vol}(T(\rho-r-\sum_{j=1}^d \jvRe \mu_j)) = C_{q,r}\rho^{d-1} + \BigO_{q,r}(\rho^{d-2} + 1).
$$
Recalling our assumptions that $|\rho| \leq 3f(h)$ and $f(h)$ may be taken sufficiently large, this proves (\ref{eIQKB}).

We introduce the further notation
\[
	A_n = \{\zeta_j \in S_1(2f(h))\::\: |\jvRe (z/h) - \jvRe \zeta_j| \in [n,n+1)\}.
\]
Our goal is to expand the product in (\ref{eIQzetaProd}) as follows:
\begin{equation}\label{eIQzetaProdExpand}
	\prod_{\zeta_j \in S_1(2f(h))} \tilde{h}^{-1}|\frac{z}{h}-\zeta_j|^{-1} = \prod_{n=0}^\infty \prod_{\zeta_j \in A_n} \tilde{h}^{-1}|\frac{z}{h} - \zeta_j|^{-1}.
\end{equation}
The set $A_n$ is contained in the union of
\[
	\{\zeta_j \in S_1(2f(h))\::\: |\rho - \jvRe \zeta_j|\leq 1/2\}
\]
for $\rho = \jvRe(z/h) + n + 1/2$ and $\rho = \jvRe(z/h) - n - 1/2$.  We may discard those $\rho$ for which $|\rho| > 3f(h)$ because we are considering only $\zeta_j$ for which $|\zeta_j| \leq 2f(h)$, making $|\rho - \jvRe \zeta_j|\leq 1/2$ impossible if $|\rho| > 3f(h)$ and $f(h)$ is large.  Using (\ref{eMuSpectrum}) and (\ref{eIQLatticeCount}), we obtain the bound
\begin{equation}\label{eIQStripBound}
	\# A_n = \BigO(f(h)^{d-1}),
\end{equation}
for all $n \geq 0$, when $f(h)$ is sufficiently large.

If $n \geq 1$, then whenever $\zeta_j \in A_n$, we have the estimate
\[
	\tilde{h}^{-1}|\frac{z}{h} - \zeta_j|^{-1} \leq \frac{f(h)}{n}.
\]
Continuing to assume $n \geq 1$, we use the upper bound (\ref{eIQStripBound}) to obtain the estimate
\[
	\prod_{\zeta_j\in A_n} \tilde{h}^{-1}|\frac{z}{h} - \zeta_j|^{-1} \leq \left\{\begin{array}{ll} \left(\frac{f(h)}{n}\right)^{Cf(h)^{d-1}}, & n < f(h) \\ 1, & n \geq f(h) \end{array}\right. .
\]
We similarly have the upper bound
\[
	\mathop{\prod_{\zeta_j\in A_0}}_{|z/h - \zeta_j| \geq 1} \tilde{h}^{-1}|\frac{z}{h} - \zeta_j|^{-1} \leq (f(h))^{Cf(h)^{d-1}}.
\]

We may combine these two estimates to bound the contribution to (\ref{eIQzetaProd}) of all $\zeta_j \in S(2f(h))$ for which $|z/h - \zeta_j| \geq 1$.  Expanding as in (\ref{eIQzetaProdExpand}), we obtain
\begin{multline}\label{eIQBigProd}
	\left(\mathop{\prod_{\zeta_j\in A_0}}_{|z/h - \zeta_j| \geq 1} \tilde{h}^{-1}|\frac{z}{h} - \zeta_j|^{-1}\right)\left(\prod_{n=1}^\infty \prod_{\zeta_j \in A_n} \tilde{h}^{-1}|\frac{z}{h} - \zeta_j|^{-1}\right)
	\\ \leq \prod_{n=1}^{\lfloor f(h) \rfloor} \left(\frac{f(h)}{n}\right)^{C f(h)^{d-1}} = \left(\frac{f(h)^{\lfloor f(h)\rfloor}}{\lfloor f(h)\rfloor !}\right)^{C f(h)^{d-1}}.
\end{multline}
From Stirling's formula, when $N$ is large,
$$
	\frac{1}{N!} \ll \frac{1}{(N/e)^N},
$$
and so, if $f(h)$ is sufficiently large,
\begin{equation}\label{eIQStirling}
	\frac{f(h)^{\lfloor f(h)\rfloor}}{\lfloor f(h) \rfloor !} \ll \left(\frac{ef(h)}{\lfloor f(h)\rfloor}\right)^{\lfloor f(h)\rfloor}.
\end{equation}
Where $f(h) \geq 2$, we have $f(h)/\lfloor f(h) \rfloor \leq 2$, and naturally $\lfloor f(h) \rfloor \leq f(h)$.  Therefore, combining (\ref{eIQBigProd}) with (\ref{eIQStirling}), we see that
$$
	\mathop{\prod_{\zeta_j \in S_1(2f(h))}}_{|z/h-\zeta_j| \geq 1} \tilde{h}^{-1}|\frac{z}{h} - \zeta_j|^{-1} \leq (2e)^{C f(h)^d}
$$
for $f(h)$ sufficiently large.  Therefore choosing $M$ sufficiently large in (\ref{eIQf}) gives one part of (\ref{eIQprod}), that
$$
	\mathop{\prod_{\zeta_j \in S_1(2f(h))}}_{|z/h-\zeta_j| \geq 1} \tilde{h}^{-1}|\frac{z}{h} - \zeta_j|^{-1} \leq \BigO(h^{-\gamma/6}).
$$

The proof is complete upon showing that
$$
	\mathop{\prod_{\zeta_j \in S_1(2f(h))}}_{|z/h-\zeta_j| < 1} \tilde{h}^{-1}|\frac{z}{h} - \zeta_j|^{-1} \leq \BigO(h^{-\gamma/6}).
$$
Applying (\ref{eIQStripBound}) to $A_0$, we see that there are $\BigO(f(h)^{d-1})$ terms in this sum.  However, we also know that quadratic operators, like the harmonic oscillator $x^2+(D_x)^2$, may in fact have eigenvalues $\lambda$ with multiplicity at least $|\lambda|^{d-1}/C$.  There is therefore no better approach, with this method, than insisting that $z/h$ be chosen a certain minimum distance from $\opnm{Spec}(q^w(x,D_x))$.

We therefore make the following assumption, with the semiclassical parameter $h$ scaled out of $q^w$:
\begin{equation}\label{eIQr}
	\opnm{dist}(\frac{z}{h}, \opnm{Spec}q^w(x,D_x)) \geq r(h).
\end{equation}
Comparing with (\ref{eIQzetaProd}), and using the counting from (\ref{eIQLatticeCount}), we find that under this assumption
$$
	\mathop{\prod_{\zeta_j \in S_1(2f(h))}}_{r(h) \leq |z/h-\zeta_j| \leq 1} \tilde{h}^{-1}|\jvRe \frac{z}{h} - \jvRe \zeta_j|^{-1} \leq \left(\frac{f(h)}{r(h)}\right)^{C_q f(h)^{d-1}}.
$$
We therefore wish to choose $r(h)$ such that
$$
	\left(\frac{f(h)}{r(h)}\right)^{C_q f(h)^{d-1}} \leq \BigO(h^{-\gamma/6}) = e^{f(h)^d/M_1},
$$
with $M_1$ possibly large depending on $M$ and $\gamma$.

Taking logarithms, we see that it is necessary to choose
$$
	\log r(h) \geq \log f(h) - \frac{f(h)}{M_1 C_q}.
$$
But $\log f(h) \ll f(h)$, and so it is sufficient to choose 
$$
	r(h) = e^{-f(h)/C_0}
$$
for any $C_0 > M_1 C_q$ when $f(h)$ is sufficiently large.  The formula (\ref{eIQSpecDist}) is a simple consequence of rescaling (\ref{eIQr}) using (\ref{RealQuadraticScaling}).

We have shown that choosing $f(h)$ as in (\ref{eIQf}) with $M$ sufficiently large but fixed, assuming distance from the spectrum of the form (\ref{eIQSpecDist}), and choosing $h_0$ sufficiently small that $f(h)$ is sufficiently large for $h \in (0,h_0]$ together establish (\ref{eIQfSmall}), (\ref{eIQexp}), and (\ref{eIQprod}).  This proves the proposition.

\end{proof}

\section{Localized quadratic resolvent estimates}\label{sLocalQuad}

In this section, we shall show that the operator $\Op{\Phi_\eps}{h}(\frakp)$, acting on $H_{\Phi_{\eps}}$, behaves very much like the unbounded elliptic quadratic differential operator
\[
\hat{q}^w(x,hD_x):L^2(\Bbb{R}^d) \rightarrow L^2(\Bbb{R}^d),
\]
provided that we localize on the FBI transform side to a neighborhood of size $|x| \leq \eps^{1/2}.$ Here the quadratic form $\hat{q}$ has been defined in (\ref{eqhat}).

Associated to $q^w(x,hD_x)$ is the corresponding quadratic operator on the FBI transform side, $\Op{\Phi_0}{h}(\frakq)$, defined by the exact Egorov relation (\ref{eEgorov}), here
\begin{equation}
\label{Q1}
T q^w(x,hD_x) = \operatorname{Op}^w_{\Phi_0, h} (\frakq)T
\end{equation}
Recall that the Weyl symbol $\frakq$ is given by $\frakq\circ \kappa_{\varphi_0}=q$. In addition to realizing $\operatorname{Op}^w_{\Phi_0, h}$ as an unbounded operator on $H_{\Phi_0}$, following \citep{H-PS}, by means of a contour deformation, we obtain the realization
$$
\Op{\Phi_q}{h}(\frakq): H_{\Phi_q}\rightarrow H_{\Phi_q}.
$$

The symbol $\frakq$ is elliptic along $\Lambda_{\Phi_q}$ in the sense that $\jvRe \frakq(x,\xi) \geq |(x,\xi)|^2/C$ for all $(x,\xi) \in \Lambda_{\Phi_q}$; this follows by comparison with $\hat{q}$ on the real side and the fact that $X \mapsto X + i\delta H_{G_q}$ is a linear isomorphism between $\Bbb{R}^{2d}$ and $\Lambda_\delta = \kap_{\varphi_0}^{-1}(\Lambda_{\Phi_q}).$  By lemma \ref{lemPhiKappaExists}, when $\delta > 0$ is sufficiently small we can find a unitary FBI transform $T_\delta$ for which
\[
	T_\delta \hat{q}^w(x,hD_x)T_\delta^* = \Op{\Phi_q}{h}(\frakq): H_{\Phi_q}\rightarrow H_{\Phi_q}.
\]

Applying proposition \ref{pImprovedQuadratic} to $\hat{q}$ and conjugating with the unitary $T_\delta$ immediately gives the following proposition.

\begin{proposition}\label{propRealSideQuadratic}
Fix $\gamma > 0.$  Define
\[
f(h) = \frac{1}{M}\left(\log \frac{1}{h}\right)^{1/d}
\]
for $M$ sufficiently large depending on $\gamma$ and $q$.  Then there exist some $C_0>0$ sufficiently large and some $h_0>0$ sufficiently small where, for any $h \in (0,h_0]$ and for any $z \in \Bbb{C}$ with $|z| \leq hf(h)$ and 
\[
	\operatorname{dist}(z, \operatorname{Spec}(\Op{\Phi_q}{h}(\frakq)) \geq he^{-f(h)/C_0},
\]
we have the resolvent estimate
\[
	||(\Op{\Phi_q}{h}(\frakq)-z)^{-1}||_{H_{\Phi_q}\rightarrow H_{\Phi_q}} \leq \BigO(h^{-1-\gamma}).
\]
\end{proposition}

\noindent\textit{Remark}. The spectrum of $\Op{\Phi_q}{h}(\frakq)$ is given by (\ref{eqnSpectrumLattice}) for the Hamilton map of the real-side quadratic form.  Because the shift of contour $\Lambda_{\Phi_0}\mapsto \Lambda_{\Phi_q}$ is between linear IR-manifolds, the corresponding Hamilton maps for $\frakq$ are obtained by conjugation with a linear canonical transformation which leaves the spectrum of the Hamilton map invariant.  We may then conclude that
$$
\operatorname{Spec}(q^w(x, hD_x)) = \operatorname{Spec}(\Op{\Phi_q}{h}(\frakq))).
$$

Following section 4 in \citep{Vi09} (which relies essentially upon section 5 in \citep{HeSjSt05}), we may next obtain localized estimates on the FBI transform side for the operator $\Op{\Phi_\eps}{h}(\frakq)$, acting on $H_{\Phi_\eps}(\Bbb{C}^d).$ We remark that, by standard formulas which follow from Fourier inversion, for $\Phi = \Phi_0, \Phi_q,$ or $\Phi_\eps$ we have the usual definition for $\Op{\Phi}{h}(\frakq)$ as a multiplication-differentiation operator:
\[
\Op{\Phi}{h}(\frakq) = \sum_{|\alpha|+|\beta|=2}(\partial^{\alpha}_x\partial^{\beta}_\xi \frakq)\left(\frac{x^{\alpha}(hD_x)^{\beta} + (hD_x)^{\beta}x^{\alpha}}{2}\right).
\]

Let $\chi_0 \in C^\infty_0(\Bbb{C}^d)$ be a cutoff function, equal to 1 in a neighborhood of 0, and $K \subset \Bbb{C}^d$ a compact neighborhood of $\operatorname{supp}\nabla\chi_0.$ The following estimates are proven in section 5 in \citep{HeSjSt05}:
\[
||(1-\Pi)\chi_0 u|| \leq \BigO(h^{1/2})||u\bar{\partial}\chi_0||
\]
and
\[
||[\Op{\Phi_q}{h}, \Pi\chi_0]u|| \leq \BigO(h)||\un_K u||
\]
for $u$ holomorphic near $\operatorname{supp}\chi_0,$ with norms in $L^2_{\Phi_q}$. Here $\Pi:L^2_{\Phi_q} \rightarrow H_{\Phi_q}$ is the
orthogonal projection onto the subspace of holomorphic functions.

We apply proposition \ref{propRealSideQuadratic} to $\Pi \chi_0 u$ when $u \in H_{\Phi_\eps}.$  When $z$ satisfies the hypotheses there, we obtain (with norms in $L^2_{\Phi_q}$)
\begin{eqnarray*}
||\chi_0 u|| &\leq& ||\Pi \chi_0 u|| + ||(1 - \Pi)\chi_0 u||
\\ &\leq& \BigO(h^{-1-\gamma})||(\Op{\Phi_q}{h}(\frakq) - z)\Pi\chi_0 u|| + \BigO(h^{1/2})||\un_K u||
\\ &\leq& \BigO(h^{-1-\gamma})||\Pi\chi_0(\Op{\Phi_q}{h}(\frakq) - z) u|| + \BigO(h^{-1-\gamma})||[\Op{\Phi_q}{h},\Pi\chi_0]u||
\\ && + \BigO(h^{1/2})||\un_K u||
\\ &\leq& \BigO(h^{-1-\gamma})||\chi_0 (\Op{\Phi_q}{h}(\frakq)-z) u|| + \BigO(h^{-\gamma})||\un_K u||.
\end{eqnarray*}

This proves the following proposition.

\begin{proposition}\label{propCutoffUnscaledQuadratic}
Fix $\gamma > 0.$  Define
\[
f(h) = \frac{1}{M}\left(\log \frac{1}{h}\right)^{1/d}
\]
for $M$ sufficiently large depending on $\gamma$ and $q$.  Let $\chi_0 \in C^\infty_0(\Bbb{C}^d)$ take values in $[0,1]$ with $\chi_0 \equiv 1$ in a neighborhood of $0 \in \Bbb{C}^d$, and let $K \subset \Bbb{C}^d$ be a fixed compact neighborhood of $\operatorname{supp}\chi_0.$ Then there exist some $C_0>0$ sufficiently large and some $h_0>0$ sufficiently small where, for any $h \in (0,h_0]$ and for any $z \in \Bbb{C}$ with $|z| \leq hf(h)$ and 
\[
	\operatorname{dist}(z, \operatorname{Spec}(\Op{\Phi_q}{h}(\frakq)) \geq he^{-f(h)/C_0},
\]
we have the localized resolvent estimate
\[
	||\chi_0 u|| \leq \BigO(h^{-1-\gamma}) ||\chi_0 (\Op{\Phi_q}{h}(\frakq) - z) u|| + \BigO(h^{-\gamma})||\un_K u||
\]
for $u \in H_{\Phi_{\eps}}$, with norms taken in $L^2_{\Phi_q}.$ Here $\un_K$ is the characteristic function of $K$.
\end{proposition}

Having localized the quadratic resolvent estimate to a fixed neighborhood of the origin in proposition \ref{propCutoffUnscaledQuadratic}, we 
must rescale to localize to a neighborhood of size $\eps^{1/2}$, on which the weights $\Phi_{\eps}$ and $\Phi_q$ agree, 
modulo higher order terms. Because both $\frakq$ and $\Phi_q$ are quadratic, we may use the rescaling (\ref{FBIChgVarsDef}) and the 
methods of section \ref{subsecRescaling}, recalling that $\tilde{h} = h/\eps$.  Applying proposition \ref{propCutoffUnscaledQuadratic} with semiclassical parameter $\tilde{h}$, we obtain, for $z$ satisfying the hypotheses of proposition \ref{propCutoffUnscaledQuadratic},
\begin{multline*}
||\chi_0(\varepsilon^{-1/2}x)u(x)||_{H_{\Phi_q,h}} = ||\chi_0 \frakU_\varepsilon u||_{H_{\Phi_q, \tilde{h}}} \\
\leq \BigO(\tilde{h}^{-1-\gamma}) ||\chi_0(\Op{\Phi_q}{\tilde{h}} (\frakq)-z)\frakU_\varepsilon u||_{H_{\Phi_q, \tilde{h}}} + \BigO(\tilde{h}^{-\gamma})||\un_K \frakU_\varepsilon u||_{H_{\Phi_q, \tilde{h}}}
\end{multline*}
Applying $\frakU_\eps^{-1}$ allows us to continue:
\begin{eqnarray*}
	&=&  \BigO(\tilde{h}^{-1-\gamma}) ||\chi_0(\varepsilon^{-1/2}\cdot) (\varepsilon^{-1}\Op{\Phi_q}{h}(\frakq)-z)u||_{H_{\Phi_q, h}} 
	\\ &&~~ + \BigO(\tilde{h}^{-\gamma}) ||\un_K(\varepsilon^{-1/2}\cdot) u||_{H_{\Phi_q,h}} \\
	&=& \BigO(h^{-1}\tilde{h}^{-\gamma}) ||\chi_0(\varepsilon^{-1/2}\cdot) (\Op{\Phi_q}{h}(\frakq)-\varepsilon z)u||_{H_{\Phi_q, h}}
	\\ && ~~ + \BigO(\tilde{h}^{-\gamma}) ||\un_K(\varepsilon^{-1/2}\cdot) u||_{H_{\Phi_q,h}}.
\end{eqnarray*}

Taking $\varepsilon = \frac{1}{C}h\log \frac{1}{h},$ we can compute that 
\[
	\tilde{h}f(\tilde{h}) = \frac{h}{\eps}\left(\frac{1}{M}\left(\log \log\frac{1}{h}\right)^{1/d} + \log \frac{1}{C}\right).
\]
Therefore, to establish $|z| \leq \tilde{h}f(\tilde{h})$, it suffices to make the assumption that 
\[
	\eps|z| \leq hF(h)
\]
for $F(h) = (1/M_0)\log\log (1/h)$ as in theorem \ref{thmMainTheorem}, with $M_0 = 2M$ and $h$ sufficiently small.  Similarly, a change of variables shows that
\[
	\operatorname{dist}(\varepsilon z, \operatorname{Spec} q^w(x, hD_x)) \geq h e^{-F(h)/C_0}
\]
suffices to establish the hypothesis $\operatorname{dist}(z, \operatorname{Spec} q^w(x, \tilde{h}D_x)) \geq \tilde{h} e^{-f(\tilde{h})/C_0}$.  The natural spectral parameter in proposition \ref{propQuadraticLocalizedFinal} below will be $\eps z.$

Since \[|(\Phi_\eps - \Phi_q)(x)| = \BigO(|x|^3),\] when both sides are localized to a region of size $\eps^{1/2},$ we have very small
difference in norm between $L^2_{\Phi_\eps}$ and $L^2_{\Phi_q}.$  More precisely, if $\operatorname{supp} v \subseteq \{|x|^2 \leq \eps\},$
then
\[
e^{-\eps^{3/2}/h} \leq \frac{||v||_{L^2_{\Phi_\eps}}}{||v||_{L^2_{\Phi_q}}} \leq e^{\eps^{3/2}/h}.
\]
Since $\eps^{3/2}/h \rightarrow 0,$ the we may replace the $L^2_{\Phi_q}$ norms in proposition \ref{propCutoffUnscaledQuadratic} with norms in $L^2_{\Phi_\eps}$ with a loss of at most a constant.

We again recall that $\Op{\Phi_q}{h}(\frakq)$ and $\Op{\Phi_\eps}{h}(\frakq)$ are identical when viewed as
quadratic forms in $(x, hD_x).$  This allows us to state our final estimate on the quadratic part of our operator.

\begin{proposition}\label{propQuadraticLocalizedFinal}
Let $\gamma > 0$ be fixed, let $\chi_0 \in C_0^\infty (\Bbb{C}^d)$ be a cutoff function taking values in $[0,1]$ with $\chi_0 \equiv 1$ near $0 \in \Bbb{C}^d$, and let
\[
F(h) := \frac{1}{C}\left(\log\log\frac{1}{h}\right)^{1/d}
\]
with $C = C(\gamma, q)$ sufficiently large.  Recall that $q$ is a quadratic form on $\Bbb{R}^{2d}$, that $\frakq = q \circ \kappa_{\varphi_0}^{-1},$ and that $\frakq$ is
elliptic along $\Lambda_{\Phi_q}.$

Assume $|z| \leq hF(h)$ and $\operatorname{dist}(z, \operatorname{Spec}(q^w(x,hD_x))) \geq h e^{-F(h)/C_0}$ for $C_0$ sufficiently large as in proposition \ref{propCutoffUnscaledQuadratic}.  Let $\eps = \eps(h) > 0$ be as in \textrm{(\ref{eqnEpsVagueDef})}, and we continue to write $\tilde{h} = h/\eps.$  We define $\Phi_\eps$ as in section \ref{ssNearWeights}.

Given these assumptions, we have the resolvent estimate
\begin{multline*}
||\chi_0(\eps^{-1/2}x)u(x)|| \leq \BigO(h^{-1}\tilde{h}^{-\gamma})||\chi_0(\eps^{-1/2}x)(\Op{\Phi_\eps}{h}(\frakq) - z)u(x)||
\\ + \BigO(\tilde{h}^{-\gamma})||\un_K(\eps^{-1/2}x) u(x)||
\end{multline*}
for any $u \in H_{\Phi_\eps}$, with norms in $L^2_{\Phi_\eps}.$
\end{proposition}

\begin{remark} As in \citep{Vi09}, we may compute the area of $\{|z| \leq hF(h)\}$ omitted by the condition 
$\operatorname{dist}(z, \operatorname{Spec}(q^w(x,hD_x))) \geq h e^{-F(h)/C_0}.$  The description of the spectrum in 
(\ref{eqnSpectrumLattice}), taken with rescaling as in (\ref{eRealRescaling}), demonstrates that
\[
\#(\operatorname{Spec}({q}^w(x,hD_x)) \cap \{|z| \leq hF(h)\}) = \BigO(F(h)^d).
\]
On the other hand, excepting circles of radius $h e^{-F(h)/C_0}$ from each of these points removes an area in $\Bbb{C}$ at most
\[
\BigO(F(h)^d)\cdot \pi (h e^{-F(h)/C_0})^2 \ll h^2 e^{-F(h)/C_0}.
\]
On the other hand, the volume of the set $\{|z| \leq hF(h)\}$ is $\pi h^2 F(h)^2$, and so
\[
\frac{\{|z|\leq hF(h)\} \cap \{\operatorname{dist}(z, \operatorname{Spec}(q^w(x,hD_x))) \geq h e^{-F(h)/C_0}\}|}{|\{|z| \leq hF(h)\}|} \ll e^{-F(h)/C_0}
\]
as $h \rightarrow 0.$  In this way, the improved estimate in proposition \ref{pImprovedQuadratic} provides that the fraction of $\{|z| \leq hF(h)\}$ to which proposition \ref{propQuadraticLocalizedFinal} does not apply is exponentially small in $F(h)$.

This will suffice to show that, at least when the spectral parameter is restricted to a region of order $h,$ the spectrum of the full operator $p^w(x, hD_x)$ approximated by the spectrum of the quadratic part $q^w(x, hD_x)$ arbitrarily closely as $h \rightarrow 0.$
\end{remark}

\section{Local estimates for full symbol}\label{sLocalFull}

Differences between the full symbol and the quadratic part at the doubly characteristic point are $\BigO(|X|^3).$  Writing $\fraka = \frakp - \frakq$ we show that, when localized, quantizations
of such symbols give small errors on $H_{\Phi_\eps}.$  While $\eps$ in the proposition is general, we will apply the proposition to $\eps$ as
in (\ref{eqnEpsVagueDef}).

\begin{proposition}\label{pLocalSmallErrors}
Let $\fraka \in S(\Lambda_{\Phi_\eps}, 1)$ be a symbol in the sense of (\ref{eFBISymbol}) obeying \[
|\fraka(X)| = \BigO(|X|^3)
\]
as $X \rightarrow 0$ along $\Lambda_{\Phi_\eps}.$ We allow $\eps > 0$ to depend on $h.$  Then, for any $u \in H_{\Phi_\eps},$
\[
||\un_{\{|x|^2\leq \eps\}}\Op{\Phi_\eps}{h}(\fraka)u||_{L^2_{\Phi_\eps}} = \BigO(\max(\eps, h)^{3/2})||u||_{H_{\Phi_\eps}}.
\]
\end{proposition}

\begin{proof} We realize $\Op{\Phi_\eps}{h}(\fraka)$ via a contour as in
(\ref{eqnFBIQuantizationIntegralCutoff}), with $
t_0$ sufficiently large
in
\begin{equation}
\label{eqnLocalFullTheta}
\theta = \frac{2}{i}(\partial_x\Phi_\eps)(\frac{x+y}{2})+it_0\overline{(x-y)},
\end{equation}
where the largeness depends on $||\nabla^2 \Phi_\eps||_{L^\infty}.$
(Recall that second derivatives of $\Phi_\eps$ are bounded independently of $h$.)  
To relate this operator to a map from $L^2(\Bbb{C}^d, dL(x))$ to itself, we first multiply an $L^2(\Bbb{C}^d, dL(x))$ function by $e^{\Phi_\eps(y)/h}$ to 
obtain a function in $H_{\Phi_\eps}$ Then we apply $\Op{\Phi_\eps}{h}(\mathfrak{a})$ to the result, and finally we multiply by $e^{-\Phi_\eps(x)/h}$ to return to $L^2(\Bbb{C}^d,dL(x)).$ The resulting integral kernel $K(x,y)$ is given by 
\[
	(2\pi h)^{-d}\fraka(\frac{x+y}{2},\theta)\jvexp \left[\frac{1}{h}\left(-\Phi_\eps(x) + i(x-y)\cdot\theta + \Phi_\eps(y)\right)\right]\psi_0(x-y).
\]

We will now prove standard estimates on the phase in $K(x,y),$ showing that
\begin{equation}
\label{eqnPhaseCancelsForSchur}
-\Phi_\eps(x) + \jvRe(i(x-y)\cdot\theta)+\Phi_\eps(y) \sim -|x-y|^2.
\end{equation}
Taylor expansion of the real-valued function $\Phi_\eps$ at $\frac{x+y}{2}$ gives
\begin{equation}\label{eqnTaylorPhiEps}
\Phi_\eps(\frac{x+y}{2} + z) = \Phi_\eps(\frac{x+y}{2}) + 2\jvRe \left(z \cdot (\partial_x\Phi_\eps)(\frac{x+y}{2})\right)  + R(x,y,z)
\end{equation}
where $|R(x,y,z)| \leq ||\nabla^2\Phi_\eps||_{L^\infty}|z|^2.$  We take the difference of (\ref{eqnTaylorPhiEps}) evaluated at $z = (y-x)/2$ and $z = (x-y)/2$ and obtain
\begin{multline*}
	\Phi_\eps(y)-\Phi_\eps(x) = 2\jvRe (\partial_x\Phi_\eps(\frac{x+y}{2})\cdot (y-x)) + \BigO(|x-y|^2)
	\\ = -\jvRe(i(x-y)\cdot\frac{2}{i}\partial_x\Phi_\eps(\frac{x+y}{2}))+\BigO(|x-y|^2)
	\\ = -\jvRe \left(i(x-y)\cdot(\theta - it_0\overline{(x-y)})\right)+\BigO(|x-y|^2).
\end{multline*}
Thus (\ref{eqnPhaseCancelsForSchur}) is established upon choosing $t_0$ sufficiently large to dominate the implicit $\Phi_\eps$-dependent constant.
$$
-\Phi_\eps(x) + 2\jvRe (\partial_x\Phi_\eps)(\frac{x+y}{2})\cdot(x-y) - t_0|x-y|^2 + \Phi_\eps(y)\sim -|x-y|^2
$$
Here we used the definition of $\theta,$ the estimate on $R$, and the assumption that $t_0$ is sufficiently large.  This proves (\ref{eqnPhaseCancelsForSchur}).

Now, from the definition of $\theta$ in (\ref{eqnLocalFullTheta}), taken with closeness of $\nabla\Phi_\eps$ to the linear $\nabla\Phi_0$ in \ref{eqnPhiDerivsClose} and the fact that $\frac{2}{i}\partial_x \Phi_\eps(0) = 0$ from the characterizations (\ref{eLambdaDeltaReal}), (\ref{eLambdaDeltaFBI}) and the fact $G_\eps(0) = 0,$ we see that \[|\theta| = \BigO(|x+y| + |x-y|).\]  We furthermore note that the estimate $|\fraka(X)| = \BigO(|X|^3)$ extends from $\Lambda_{\Phi_\eps}$ to all of $\Bbb{C}^{2d}$ when $\fraka$ is extended almost holomorphically off $\Lambda_{\Phi_\eps}$.  Using in addition
that $|x+y| \leq 2|x| + |x-y|,$ we get that \[\fraka(\frac{x+y}{2}, \theta) = \BigO(|x|^3 + |x-y|^3).\] Recall as well that $dy\wedge
d\theta = \BigO(1)dL(y).$ This gives a kernel for the integral operator on $L^2(\Bbb{C}^n, dL(x))$ bounded by
\[\BigO(1)h^{-d}e^{-\frac{1}{Ch}|x-y|^2}\left(|x|^3 + |x-y|^3\right)\un_{\{|x|^2\leq \eps\}}(x).\]  The proposition then follows from Schur's
test, since \[\iint |x|^3\un_{\{|x|^2\leq \eps\}}(x)h^{-d}e^{-\frac{1}{Ch}|x-y|^2}\,\left(dL(x)~\mathrm{or}~dL(y)\right) =
\BigO(\eps^{3/2})\] from $|x| \leq \eps^{1/2}$ on the support of the integrand and since \[ \iint \un_{\{|x|^2\leq \eps\}}(x)h^{3/2}h^{-d}
\left(\frac{|x-y|}{h^{1/2}}\right)^3 e^{-\frac{1}{Ch}|x-y|^2}\,\left(dL(x)~\mathrm{or}~dL(y)\right) = \BigO(h^{3/2})\] from a change of variables.

\end{proof}

\section{Estimates for exterior region}\label{sExterior}
Here we shall establish resolvent type estimates localized to the region outside a tiny $h$-dependent neighborhood of the
doubly characteristic point. As before, we shall consider the IR-manifold
$$
\Lambda_{\Phi_{\eps}}=\kappa_{\varphi_0}(\Lambda_{\delta,\eps})=
\Big\{\Big(x,\frac{2}{i}\frac{\partial \Phi_{\eps}}{\partial x}(x)\Big) : x \in \Bbb{C}^d\Big\},\
$$
associated to the weight $G_{\eps}$, and we recall that the small parameter $\eps$ is taken as in (\ref{eqnEpsVagueDef}) equal to
$$
\eps=\frac{1}{C} h\log \frac{1}{h},
$$
where $C > 0$ may be taken large but will be fixed in the proof of the theorem.
We shall be concerned with studying the region on the FBI-transform side of the IR-manifold $\Lambda_{\Phi_{\eps}}$ where
\begin{equation}
\label{ext1}
|x| \geq \eps^{1/2}.
\end{equation}
Note also that we have, in this region,
\begin{equation}
\label{ext2}
\xi = -\jvIm x +\BigO(\delta \eps^{1/2})
\end{equation}
for $(x,\xi)\in \Lambda_{\Phi_{\eps}}$. When working in the unbounded region (\ref{ext1}), we recall from Proposition 3.2 that we have
\begin{equation}
\label{ext3}
\jvRe \frakp \Big(x,\frac{2}{i}\frac{\partial \Phi_{\eps}(x)}{\partial x}\Big)\Big) \geq \frac{\delta \eps}{\tilde{C}},
\end{equation}
when $|x| \geq \eps^{1/2}$. It is therefore convenient to consider the new rescaled variables
\begin{equation}
\label{ext4}
x = \eps^{1/2} \widetilde{x}.
\end{equation}
Working in the rescaled variables, we shall show how to obtain the following result.

\begin{proposition}\label{propExteriorRegion}
Let $\chi \in C^\infty(\Bbb{C}^d)$ be fixed, taking values in $[0, 1]$, equal to zero in a neighborhood of $0 \in \Bbb{C}^d,$ and equal to $1$ off
a compact set.  Assume that $\frakp$ and $z$ continue to satisfy the hypotheses of theorem \textrm{\ref{thmMainTheorem}}.  
Then, for $\eps$ as in \textrm{(\ref{eqnEpsVagueDef})} and $u \in H_{\Phi_\eps},$ we have
\begin{multline}
\label{eqnExteriorEstimate}
\int_{\Bbb{C}^d} \chi(\eps^{-1/2}x)|u|^2 e^{-\frac{2}{h}\Phi_\eps}\,dL(x)
\\ \leq \BigO(\eps^{-1})||(\Op{\Phi_\eps}{h}(\frakp)-z)u||\:||u|| + \BigO(\tilde{h})||u||^2,
\end{multline}
with norms taken in $H_{\Phi_\eps},$ and $\tilde{h} = h/\eps.$
\end{proposition}

\begin{proof}
We use the rescaling formulas in section \ref{subsecRescaling}, recalling that $\tilde{\Phi}_\eps(x) = \eps^{-1}\Phi_\eps(\eps^{1/2}x)$.  We obtain
\[
\langle \chi(\eps^{-1/2}x)(\Op{\Phi_\eps}{h}(\frakp)-z)u,u\rangle_{H_{\Phi_\eps, h}} = \langle \chi(x)\Op{\tilde{\Phi}_\eps}{\tilde{h}}(\frakp_\eps - z)\frakU_\eps u, \frakU_\eps u\rangle_{H_{\tilde{\Phi}_\eps,\tilde{h}}}.
\]
Here $\frakp_\eps(x,\xi) = \frakp(\eps^{1/2}x,\eps^{1/2}\xi),$ and multiplying by $\eps^{-1}$ gives us a symbol $\eps^{-1}(\frakp_\eps - z)$ such that
$$
|\partial^\alpha_{x,\xi}\left(\eps^{-1}(\frakp_\eps - z)\right)| = \BigO_\alpha(1),
$$
uniformly with respect to $\eps > 0$, when $|\alpha| \geq 2$. It follows from (\ref{ext3}) that, along the rescaled manifold given by
$$
\tilde{\xi}(\tilde{x}) = \frac{2}{i}(\partial_{\tilde{x}}\tilde{\Phi}_\eps)(\tilde{x})
$$
and restricted to the region where $\tilde{x}\in \operatorname{supp}\chi$, the real part of $\eps^{-1}\frakp_\eps$ is uniformly
bounded from below by a fixed positive constant. Furthermore, restricting the attention to the same manifold and recalling that $(0,0)$ is a doubly characteristic point for $\frakp$, we have
$$
\frac{1}{\eps}\frakp_\eps = \BigO(1+|\tilde{x}|^2)
$$
uniformly in $\eps > 0$.
Since $|z| \leq hF(h) \ll \eps,$ the uniform lower bound for $\eps^{-1} \frakp_\eps$ implies a uniform lower bound for $\eps^{-1}(\frakp_\eps - z)$.

Applying the quantization-multiplication formula (\ref{QM}) in the rescaled variables and taking real parts then gives
\begin{multline*}
\eps^{-1}\jvRe \langle \chi(x)(\Op{\tilde{\Phi}_\eps}{h}(\frakp_\eps) - z)\frakU_\eps u, \frakU_\eps u\rangle_{H_{\tilde{\Phi}_\eps,\tilde{h}}}
\\ \geq \int \chi(x)\eps^{-1}\jvRe(\frakp_\eps(x, \tilde{\xi}(x))-z)|(\frakU_\eps u)(x)|^2 e^{-2\tilde{\Phi}_\eps(x)/\tilde{h}}\,dL(x)
\\ - \BigO(\tilde{h})||\frakU_\eps u||^2_{H_{\tilde{\Phi}_\eps,\tilde{h}}}.
\end{multline*}
Here, as above, $\tilde{\xi}(x) = \frac{2}{i}(\partial_x\tilde{\Phi}_\eps)(x).$

Using the lower bound for $\eps^{-1}\jvRe(\frakp_\eps(x,\tilde{\xi}(x))-z),$ changing variables via the rescaling $\frakU_\eps^{-1}$, and
using Cauchy-Schwarz on the inner product gives the conclusion, (\ref{eqnExteriorEstimate}).

\end{proof}

\section{Proof of theorem}\label{sProof}

Following section 6 in \citep{Vi09} or section 7 in \citep{H-PS}, we glue together exterior and interior estimates to create a resolvent estimate and prove theorem \ref{thmMainTheorem}.

Let $u$ be an arbitrary element of $H_{\Phi_\eps},$ the weighted space described in section \ref{ssNearWeights}, with $\eps =
\frac{1}{C}h\log\frac{1}{h}$ and $\tilde{h} = h/\eps$.  We assume that $\frakp = p \circ \kap_{\varphi_0}^{-1}$ is an almost analytic extension of
$\frakp$ off $\Lambda_{\Phi_0}$ and that $\frakq = q \circ \kap_{\varphi_0}^{-1}$ is the corresponding holomorphic quadratic approximation near the doubly characteristic point $(0,0)$.
We also continue to assume that the spectral parameter $z$ satisfies the assumptions of theorem \ref{thmMainTheorem}: that $|z| \leq hF(h)$ for 
\[
	F(h) = \frac{1}{C_0} \left(\log \log \frac{1}{h}\right)^{1/d}
\]
and
\[
	\operatorname{dist}(z, \operatorname{Spec}(q^w(x,hD_x))) \geq h e^{-F(h)/C_1}
\]
for $C_0, C_1$ sufficiently large.

Let $\chi_0(x) \in C^\infty_0(\Bbb{C}^d)$ be a cutoff function taking values in $[0,1]$, with $\chi(x)$ equal to $1$ for $|x| \leq 1/2$ and equal to zero for $|x| \geq 1.$  As in proposition
\ref{propQuadraticLocalizedFinal}, let $K$ be a compact neighborhood of $\operatorname{supp} \chi_0$ avoiding $0 \in \Bbb{C}$; we may, for example, say $K = \{1/3 \leq |x| \leq 3/2\}.$  Let
$\chi_1(x) \in C^\infty(\Bbb{C}^d)$  be a cutoff function taking values in $[0,1]$ and localizing to a neighborhood of infinity containing $K$, meaning for instance that $\chi_1(x)$ equals zero for $|x| \leq 1/4$ and equals $1$ for $|x| \geq 1/3.$

For brevity of notation, we let $\frakp^w$ denote $\Op{\Phi_\eps}{h}(\frakp)$ and $\frakq^w$ denote $\Op{\Phi_\eps}{h}(\frakq).$  To denote rescaled cutoff functions, we let $\tilde{\chi}_0(x) = \chi_0(\eps^{-1/2}x)$, and we define $\tilde{\chi}_1$ and $\tilde{\un}_K$ analogously.  Unless otherwise stated, we assume norms are in $L^2_{\Phi_\eps}.$

Therefore, using Proposition \ref{propQuadraticLocalizedFinal}, Proposition \ref{pLocalSmallErrors} applied to $\fraka = \frakp -
\frakq$, and the observation that $\un_K \leq \chi_1,$ we obtain
\begin{eqnarray*}
||u|| &\leq& ||\tilde{\chi}_0 u|| + ||\tilde{\chi}_1 u||
\\ &\leq& \BigO(h^{-1}\tilde{h}^{-\gamma})||\tilde{\chi}_0(\frakq^w-z)u|| + \BigO(\tilde{h}^{-\gamma})||\tilde{\un}_K u|| + ||\tilde{\chi}_1 u||
\\ &\leq& \BigO(h^{-1}\tilde{h}^{-\gamma})\left(||\tilde{\chi}_0(\frakp^w-z)u|| + ||\tilde{\chi}_0(\frakq^w - \frakp^w)u||\right) + \BigO(\tilde{h}^{-\gamma})||\tilde{\chi}_1 u||
\\ &\leq& \BigO(h^{-1}\tilde{h}^{-\gamma})||(\frakp^w-z)u|| + \BigO(h^{-1}\tilde{h}^{-\gamma}\eps^{3/2})||u|| + \BigO(\tilde{h}^{-\gamma})||\tilde{\chi}_1 u||.
\end{eqnarray*}
Since $\eps$ is only logarithmically larger than $h,$ the rescaled parameter $\tilde{h}$ is only logarithmic in $h$ while $\eps^{3/2}$ is nearly as small as $h^{3/2}.$  Therefore the second term only logarithmically larger than $\BigO(h^{1/2})||u||$ and can be absorbed into the left-hand side when $h$ is small with no difficulty.

To bound $||\tilde{\chi}_1 u||,$ we apply proposition \ref{propExteriorRegion} with cutoff function $\tilde{\chi} = \tilde{\chi}_1^2,$ and note 
that the left-hand side of (\ref{eqnExteriorEstimate}) is then precisely $||\tilde{\chi}_1 u||^2.$  Thus
\[
\BigO(\tilde{h}^{-\gamma})||\tilde{\chi}_1 u|| 
\leq \BigO(\tilde{h}^{-\gamma})\left[\BigO(\eps^{-1})||(\frakp^w - z)||\:||u|| + \BigO(\tilde{h})||u||^2\right]^{1/2}.
\]
Distributing the square root at the loss of a constant and using the Cauchy-Schwarz inequality gives
\[
\BigO(\tilde{h}^{-\gamma})||\tilde{\chi}_1 u|| \leq \BigO(\tilde{h}^{-\gamma}\eps^{-1})||(\frakp^w - z)u|| + (\frac{1}{2} + \BigO(\tilde{h}^{1/2 - \gamma}))||u||.
\]
Since $\eps \gg h,$ we have $\tilde{h}^{-\gamma}\eps^{-1} \ll h^{-1}\tilde{h}^{-\gamma},$ so we may absorb the first term into $\BigO(h^{-1}\tilde{h}^{-\gamma})||(\frakp^w-z)u||$ in the previous estimate for the entirety of $||u||.$  So long as $\gamma < 1/2,$ the coefficient of $||u||$ immediately above is less than $1$ for $h$ sufficiently small.  Therefore, with the restriction $\gamma < 1/2,$ we absorb the term involving $||u||$ into the left-hand side of the estimate for the entirety of $||u||$ and obtain
\[
||u|| \leq \BigO(h^{-1}\tilde{h}^{-\gamma})||(\frakp^w-z)u||.
\]

Using (\ref{eqnNormsClose}), we replace norms in $H_{\Phi_\eps}$ with norms in $H_{\Phi_0}$ at the price of a multiplicative factor $e^{\pm \BigO(1)\eps/h}.$  We may increase the constant
in the definition (\ref{eqnEpsVagueDef}) of $\eps$, depending only on $\rho$ and the $h$-independent bound (\ref{eqnPhisClose}), to obtain $e^{\BigO(1)\eps/h} = \BigO(h^{\rho/3})$.  This suffices to show that
\[
||u||_{H_{\Phi_0}} \leq \BigO(h^{-1-\rho})||(\frakp^w-z)u||_{H_{\Phi_0}}.
\]
Recalling (\ref{error}), we have that $\frakp^w = \Op{\Phi_\eps}{h}(\frakp)$ differs from $\Op{\Phi_0}{h}(\frakp)$ by $\BigO(h^\infty):L^2_{\Phi_0}\rightarrow
L^2_{\Phi_0}.$ This establishes the real-side resolvent estimate in theorem \ref{thmMainTheorem} after conjugation with the standard unitary FBI transform
$T_{\varphi_0}:L^2(\Bbb{R}^d)\rightarrow H_{\Phi_0}(\Bbb{C}^d)$ as a consequence of the exact Egorov theorem (\ref{eEgorov}).

The a priori resolvent estimate implies existence of the resolvent by virtue of
$$
\{p^w(x,hD_x) - z\::\: z \in \operatorname{neigh}(0; \Bbb{C})\}
$$
being a analytic family of Fredholm operators of index 0 for $h$ sufficiently small.  The proof of theorem \ref{thmMainTheorem} is now complete.

\bibliography{MicrolocalBibliography}

\end{document}